\documentclass[a4paper,fleqn]{cas-sc}

\usepackage[authoryear,longnamesfirst]{natbib}

\usepackage{hyperref}
\usepackage{amsmath}
\usepackage{amssymb}
\usepackage{amsthm}
\usepackage{graphicx}
\usepackage{graphics}
\usepackage{tikz-cd}
\usepackage{physics}
\usepackage{tikz}
\usetikzlibrary{arrows.meta}
\usetikzlibrary{decorations.pathmorphing}
\usetikzlibrary{arrows}

\def\tsc#1{\csdef{#1}{\textsc{\lowercase{#1}}\xspace}}
\tsc{WGM}
\tsc{QE}

\newtheorem{theorem}{Theorem}
\newtheorem{lemma}[theorem]{Lemma}
\newdefinition{rmk}{Remark}
\newproof{pf}{Proof}
\newproof{pot}{Proof of Theorem \ref{thm}}

\newdefinition{definition}{Definition}
\newtheorem{corollary}[theorem]{Corollary}
\newtheorem{proposition}[theorem]{Proposition}
\newdefinition{ex}{Example}
\newdefinition{conj}{Conjecture}

\newcommand{\kk}{k\langle \langle x_0, x_1 \rangle \rangle}
\DeclareMathOperator{\conc}{conc}
\newcommand{\Z}{\mathbb{Z}}
\DeclareMathOperator{\Li}{Li}
\DeclareMathOperator{\dep}{dp}
\DeclareMathOperator{\Sh}{Sh}
\DeclareMathOperator{\dpt}{dp}
\DeclareMathOperator{\pr}{pr}
\newcommand{\N}{\mathbb{N}}
\DeclareMathOperator{\ad}{ad}
\DeclareMathOperator{\id}{id}
\newcommand{\g}{\mathfrak{g}}
\DeclareMathOperator{\kv}{kv}
\DeclareMathOperator{\wt}{wt}

\newcommand{\grt}{\mathfrak{grt}}
\newcommand{\rc}{\mathfrak{rc}}
\newcommand{\dmr}{\mathfrak{dmr}}
\newcommand{\krv}{\mathfrak{krv}}

\newcommand{\fr}{\mathfrak{fr}}

\newcommand{\frone}{\mathfrak{fr}^{>1}_k(x_0,x_1)}
\newcommand{\ldb}{\{\!\!\{}
\newcommand{\rdb}{\}\!\!\}}

\DeclareMathOperator{\sk}{Skew}

\ExplSyntaxOn
\bool_gset_true:N \g_stm_nologo_bool
\ExplSyntaxOff

\begin{document}
\let\WriteBookmarks\relax
\def\floatpagepagefraction{1}
\def\textpagefraction{.001}

\shorttitle{}    

\shortauthors{M. Howarth, M. Ren}  

\title [mode = title]{Polylogarithmic characterizations of the reduced coaction Lie algebra and the double shuffle Lie algebra, and their relation to the Kashiwara--Vergne Lie algebra}  



%

\author[1]{Megan Howarth}

\cormark[1]


\ead{megan.howarth@unige.ch}

\ead[url]{}


\affiliation[1]{organization={Université de Genève},
            addressline={Rue du Conseil-Général 7-9}, 
            city={Genève},
          citysep={}, 
            postcode={1205},
            country={Switzerland}}

\author[2]{Muze Ren}


\ead{muze.ren@unige.ch}

\ead[url]{}

\credit{}

\affiliation[2]{organization={Institut de Recherche Mathématique Avancée, UMR 7501, Université de Strasbourg et CNRS},
            addressline={Rue René Descartes 7}, 
            city={Strasbourg},
            postcode={67000}, 
            country={France}}

\cortext[1]{Corresponding author}



\begin{abstract}
The reduced coaction Lie algebra $\rc_0$ is defined by a
skew-symmetric condition together with an algebraic equation involving
the reduced coaction, a non-cyclic refinement of the necklace
cobracket. The double shuffle Lie algebra $\dmr_0$ encodes
the regularized double shuffle relations in the study of formal multiple zeta
values, while the Kashiwara--Vergne Lie algebra $\krv_2$ arose in the study of the Kashiwara--Vergne conjecture.

Our main result is the description of $\dmr_0$ and $\rc_0$ by vanishing conditions of different strengths on the same object: the defect of Drinfeld's pentagon equation; $\dmr_0$ is characterised by the vanishing of the associated two-variable multiple polylogarithms, and $\rc_0$ by the vanishing of those whose second index has depth $1$. 

We also establish explicit relations among these three Lie algebras. First, the skew-symmetric elements of $\mathfrak{dmr}_0$ embed into
$\rc_0$. Second, under an additional condition, we construct an injective Lie algebra morphism from $\rc_0$ to $\krv_2$.
This exhibits $\rc_0$ as a natural intermediate structure between
double shuffle and Kashiwara--Vergne theory.
\end{abstract}



\maketitle









\tableofcontents
\section{Introduction and main results}
\subsection{The Grothendieck-Teichm\"uller, double shuffle and Kashiwara--Vergne Lie algebras}
\label{subsec:three_Lie_algebra}

The Grothendieck-Teichmüller, double shuffle and Kashiwara--Vergne Lie algebras were introduced in different contexts and are each of independent interest; we start by reviewing them. \\
Multiple zeta values (MZVs) satisfy two families of product relations: the
iterated-integral representation gives the \textit{shuffle} relations, whereas
the series representation gives the \textit{stuffle} 
relations. After a suitable regularization, the combination of these
two families is known as the \emph{regularized double shuffle
relations}. Racinet \cite{Racinet2002} introduced the \textit{double shuffle Lie algebra} $\mathfrak{dmr}_0$
as the Lie algebra of formal solutions to these relations, equipped with
the Ihara bracket. It plays a central role in the algebraic study of formal MZVs.

The \textit{Grothendieck-Teichmüller Lie algebra} $\mathfrak{grt}_1$ introduced by Drinfeld in \cite{Drinfeld1991} plays a pivotal role in many branches of mathematics. In particular, it describes the infinitesimal symmetries of Drinfeld associators. It is characterised by three defining equations: skew symmetry, infinitesimal hexagon and pentagon; see equations (5.17), (5.18) and (5.20) in \cite{Drinfeld1991}, respectively. It turns out that the pentagon equation is the only necessary equation, as Furusho proved that the pentagon equation implies the two others \cite{furusho_pentagon_2010}. The study of the pentagon equation, or its \textit{defect}, will play a central role in this work. In \cite{Furusho2011}, Furusho also showed that Drinfeld associators satisfy the regularized double shuffle relations; at the infinitesimal level, this yields an injective Lie algebra morphism $\grt_1 \hookrightarrow \dmr_0$. The polylogarithmic interpretation of this result is particularly important for us: suitable one- and two-variable multiple polylogarithms detect the double shuffle relations through the defect of the pentagon equation.

The third Lie-theoretic structure arises from the Kashiwara--Vergne theory. Originally concerning the Duflo isomorphism, Alekseev and Torossian formulated in \cite{Alekseev2012} its infinitesimal version in terms of the \textit{Kashiwara--Vergne Lie algebra} $\mathfrak{krv}_2$,
whose elements are special tangential derivations satisfying a divergence condition. The recent works of Schneps \cite{Schneps2025}, Enriquez and Furusho \cite{EF4} and Kawamura \cite{kawamura2025} independently relate the double shuffle Lie algebra $\dmr_0$ to $\krv_2$. 
It is an open question whether these three Lie algebras are isomorphic.

\subsection{Goldman-Turaev theory and the reduced coaction Lie algebra} 
\label{subsec:coaction-background}
For the surface $\Sigma_{0,n+1}$ of genus zero with $n+1$ boundary components, let $\g(\Sigma)$ denote the vector space spanned by free homotopy classes of loops. Equipped with the Goldman bracket and Turaev cobracket, it yields the Lie bialgebra $\left( \g(\Sigma),[\cdot,\cdot]_{\mathrm{Goldman}}, \delta^f_{\mathrm{Turaev}} \right)$, where the augmentation ideal filtration on $\pi_1(\Sigma_{0,n+1})$ moreover induces a filtration on $\g(\Sigma)$. The associated graded is the necklace Lie bialgebra $\mathrm{gr}(\g(\Sigma), [-,-]_{\mathrm{necklace}}, \delta_{\mathrm{necklace}})$ associated to the star-shaped quiver.

Both the topological and algebraic cobrackets admit refinements called
\emph{reduced coactions}. The essential point is that the cobracket
records two cyclic outputs obtained by cutting a self-intersection,
whereas the coaction retains one output as a based path or,
algebraically, as a non-cyclic word. Thus, the reduced coaction carries
more information than the corresponding cobracket. One then recovers the latter after applying the
appropriate cyclic projection and antisymmetrization, and recovers the corresponding bracket through the product formula; see Section 3.2 of \cite{ANR} for more details. Therefore, the reduced coaction encodes the information of both the bracket and of the cobracket. \\

We briefly recall the precise meaning of \textit{coaction} and \textit{reduced coaction}; for additional background, see the notions of \textit{coaction} in \cite{Kawazumi2015} and \textit{quasi derivation} in \cite{Massuyeau2018}. Let $A=(A, \cdot, \Delta, 1, S)$ be an involutive Hopf algebra, i.e. $S^2=id$, a map $\bar{\mu}: A \to |A| \otimes A$ is called a \emph{coaction} map with respect to the double bracket $\ldb -, - \rdb$ (in the sense of  Van den Bergh \cite{double_poisson}) if it satisfies the product rule
\begin{equation}   \label{eq:mu_mult}
	\bar{\mu}(ab)=\bar{\mu}(a)(1\otimes b)+(1\otimes a)\bar{\mu}(b)+(|-|\otimes 1)\ldb a,b\rdb.
\end{equation}
We define the \emph{reduced coaction map} $\mu$ associated to $\bar{\mu}$ by 
$\mu=(\varepsilon \otimes 1)\bar{\mu}: A \to A$. Note that for the completed free associative algebra $A=k\langle\!\langle x_0,x_1\rangle\!\rangle$,
the relevant double bracket is the KKS double bracket.

In \cite{GTgenus0}, Alekseev, Kawazumi, Kuno and Naef related the genus-zero formality
problem for the Goldman-Turaev Lie bialgebra to the
Kashiwara--Vergne problem. Their work also explains the
close relation between the necklace cobracket and the divergence of a
tangential derivation. At the refined infinitesimal level, the algebraic reduced
coaction may be regarded as a non-cyclic lift of the data that yields
the divergence after passage to cyclic words; see Section \ref{sec:Kashiwara--Vergne}. \\

Motivated by these constructions and the importance of the reduced coaction, the second author introduced the reduced coaction Lie algebra $\overline{\mathfrak{rc}_0}$ in \cite{Ren2025}. It consists of suitably normalized skew-symmetric Lie series satisfying an equation expressed in terms
of the reduced coaction map associated to the KKS double bracket. It is therefore natural to ask
whether the reduced coaction Lie algebra
$\overline{\rc}_0$, or the equivalent form $\rc_0$ we consider, provides an 
intermediate structure
between the double shuffle theory and the Goldman-Turaev approach to the
Kashiwara--Vergne problem.

\subsection{Motivation and contributions of the paper}

The preceding developments provide two different approaches to the
Kashiwara--Vergne theory. The first is arithmetic and
associator-theoretic: $\grt_1 \hookrightarrow \dmr_0 \hookrightarrow \krv_2.$ 
The second originates in the topology of surfaces and passes through
the Goldman-Turaev Lie bialgebra, the necklace Lie bialgebra,
reduced coactions and divergence. The reduced coaction Lie algebra $\rc_0$ naturally belongs to the second, and in this paper, we show how it relates to the first. \\

The starting point is that the aforementioned Lie algebras $\dmr_0$ and $\rc_0$ can be characterized by vanishing conditions on the same object: the defect of Drinfeld's pentagon equation, evaluated against two-variable multiple polylogarithms. We show that $\rc_0$ is also characterised by such a condition: the vanishing in the skew-symmetric setting of two-variable multiple polylogarithms whose second index has depth one.
The conditions we impose are of successively decreasing strength, which places the Lie algebras into a hierarchy, and in this sense, $\rc_0$ is an intermediate structure between $\dmr_0$ and the Goldman-Turaev approach to the Kashiwara--Vergne problem, answering the question of the previous Subsection \ref{subsec:coaction-background}. 

Furthermore, the reduced coaction equation is a non-cyclic refinement of the Kashiwara--Vergne divergence condition. After symmetrizing and passing to cyclic words, the reduced coaction $\mu$ recovers the divergence of the associated special derivation. 
It follows that, subject to an additional equation (Vep), an element of $\rc_0$ determines an element of $\krv_2$, yielding an injective Lie algebra morphism between the two spaces.

Thus, the role of the reduced coaction is twofold: it isolates a natural depth-one portion of the
double shuffle constraints and, at the same time, provides a lift of
the divergence condition appearing in the Kashiwara--Vergne theory. \\

The main tools throughout are polylogarithmic computations in the spirit of
Furusho \cite{furusho_pentagon_2010}, Fox derivatives and Fox pairings, the KKS double
bracket, and cyclic-word calculus.

\subsection{Main results}
We fix $k$ to be a field of characteristic zero and consider the Hopf algebra of formal noncommutative power series 
\begin{equation}
\label{eq:def-hopf-algebra-A}
    A := 
    \big( \kk, \conc, \Delta_\shuffle, \varepsilon,S \big),
\end{equation}
where $\conc$ is the concatenation product and for $i=0,1$, 
\begin{equation*}
    \begin{split}
        & \Delta_\shuffle(x_i) = x_i \otimes 1 + 1 \otimes x_i, \quad \varepsilon(x_i) = 0, \quad S(x_i)=-x_i. 
    \end{split}
\end{equation*}

Every element $\varphi\in A$ 
is uniquely expressed as
\begin{equation}\label{eq:decomposition}
\varphi =\varepsilon(\varphi)+ x_0 d^R_0(\varphi) + x_1 d^R_1(\varphi)
    =\varepsilon(\varphi)+ d^{L}_0(\varphi) x_0 + d^{L}_1 (\varphi)x_1,
\end{equation}
where 
$d_i^R(\varphi)$ and $d_i^L(\varphi)$ denote the remaining parts of $\varphi$ after removing respectively the initial or final letter $x_i$, for $i \in \{0,1\}$. We denote by $c_{w}(\varphi)$ the coefficient of the word $w$ in $\varphi$. \\

Let $\mathfrak{fr}_k(x_0,x_1)$ be the free Lie algebra in $k\langle\langle x_0,x_1\rangle\rangle$ generated by $x_0,x_1$ and $\mathfrak{fr}^{>1}_k(x_0,x_1)$ be the elements such that $\psi(x_0,0)=\psi(0,x_1)=0$. We introduce the following two sub-vector spaces for later use:
\begin{align}
 \mathrm{Skew} & :=\{\psi\in \mathfrak{fr}^{>1}_k(x_0,x_1)\mid \psi(x_0,x_1)=-\psi(x_1,x_0)\}, \label{def:Skew} \\
 \mathrm{Vep} & :=\{\psi\in \mathfrak{fr}^{>1}_k(x_0,x_1)\mid [x_1,\psi(-x_0-x_1,x_1)]+[x_0,\psi(-x_0-x_1,x_0)]=0\}. \label{def:Vep}
\end{align}

\begin{rmk} \label{rem:skew}
The skew symmetry condition considered here is neither the symmetric condition of \cite{Alekseev2012}, nor the involution $\Theta$ of \cite{EF4}; it is however directly related to the $S_3$-symmetry appearing in \cite{Drinfeld1991} and \cite{Ihara1992}. The Vep condition corresponds to the equation $(5.19)$ in \cite{Drinfeld1991}, which is a variant of the special derivation property discussed in Section \ref{sec:Kashiwara--Vergne}.
\end{rmk}

In Section \ref{section:notation_polylogarithms}, we introduce the basics of the formal two-variable polylogarithms $l^{y,x}_{\bf{a},\bf{b}}$ and their variants. The main use of these polylogarithms is that, when evaluated in the defect $\alpha$ of a variant of Drinfeld's pentagon equation \cite[Equation (5.20)] {Drinfeld1991}, they detect structural information about related Lie algebras. Here, we consider it in the form 
\begin{equation*}
\alpha\left(\psi(x_0,x_1)\right):=\psi_{451}+\psi_{123}-\psi_{432}-\psi_{215}-\psi_{543},
\end{equation*} 
where $\psi_{ijk}:=\psi(x_{ij}, x_{jk})$, valued in $\mathfrak{p}_5$ (see Definition \ref{def:pure-shere-braid}). The relation between the pentagon equation and two-variable polylogarithms was first considered by Furusho in \cite{Furusho2011} for the double shuffle Lie algebra $\dmr_0$. We continue the study of $\dmr_0$ in Section \ref{sec:poly-descr-double-shuffle}. 

\begin{definition} \label{def:double-shuffle}
The \textit{double shuffle Lie algebra} $\mathfrak{dmr}_0$ is defined to be the set of formal Lie series $\psi\in \mathfrak{fr}_k(x_0,x_1)$ satisfying  
\begin{equation*} \label{eq:def-dmr0}
	c_{x_0}(\psi)=c_{x_1}(\psi)=c_{x_0x_1}(\psi)=0 \quad \text{ and } \quad \Delta_{*}(\psi_{*})=1\otimes \psi_{*}+\psi_{*}\otimes 1,
\end{equation*}
where $\psi_{*}=\psi_{\text{corr}}+\pi_{Y}(\psi)$, with $\psi_{\text{corr}}:=\sum_{n=1}^{\infty}\frac{(-1)^n}{n}c_{x^{n-1}_0x_1}(\psi)y^n_1$ and
$\pi_Y$ the $k$-linear map defined by
\begin{equation*}
    \begin{split}
        k\langle\langle x_0,x_1\rangle\rangle & \to k\langle\langle y_1,\dots,y_n,\dots\rangle\rangle \\
        \text{words ending in } x_0 & \mapsto 0\\
        x^{n_m-1}_0x_1\dots x^{n_1-1}_0x_1 & \mapsto (-1)^m y_{n_m}\dots y_{n_1}.
    \end{split}
\end{equation*}

The coproduct $\Delta_{*}$ on $k\langle\langle y_1,y_2,\dots \rangle\rangle$ is defined to be $\Delta_{*}y_n=\sum_{i=0}^n y_i\otimes y_{n-i}$, with $y_0=1$.
\end{definition}

We recall that one of the fundamental results proved by Racinet \cite{Racinet2002} about $\mathfrak{dmr}_0$ is that it is a Lie algebra with respect to the \textit{Ihara bracket}
	\begin{equation}\label{eq:Ihara_bracket}
		\{\psi_1,\psi_2\}=d_{\psi_2}(\psi_1)-d_{\psi_1}(\psi_2)-[\psi_1,\psi_2],
	\end{equation}
	where $d_{\psi}$ is the derivation of $\psi$ 
    given by $d_{\psi}(x_0)=0$ and $d_{\psi}(x_1)=[x_1,\psi]$. \\

Our first result is a new polylogarithmic characterization of $\mathfrak{dmr}_0$ in terms of the defect $\alpha$ 
of the pentagon equation  

\renewcommand{\thetheorem}{\Alph{theorem}}
\setcounter{theorem}{0}
\begin{theorem}[Theorem \ref{theorem:dmr_defect}] \label{theorem:dmr_defect-B}
	Let $\psi\in \mathfrak{fr}_k(x_0,x_1)$ be such that $c_{x_0}(\psi)=c_{x_1}(\psi)= c_{x_0 x_1}(\psi)=0$, then the following two conditions are equivalent:
	\begin{enumerate}[(i)] 
		\item $\psi\in \mathfrak{dmr}_0$;
		\item $l^{y,x}_{\bf{a},\bf{b}}(\alpha)=0,$ for ${\bf{a},\bf{b}}\ne (1,\dots,1),(1,\dots,1)$.
	\end{enumerate}
\end{theorem}

In Section \ref{sec:poly-descr-reduced-coaction}, we undertake a similar study for the polylogarithmic description of the reduced coaction Lie algebra $\mathfrak{rc}_0$, which is formally defined as follows.

\begin{definition}[The reduced coaction map]
The (algebraic) \textit{reduced coaction} $\mu$ is a linear map from $k\langle\langle x_0,x_1\rangle\rangle$ to $k\langle\langle x_0,x_1\rangle\rangle$ defined by
\begin{equation}\label{eq:reduced_coaction}
           \mu(x_0)=\mu(x_1)=0, \quad
           \mu(k_1\,k_2\dots\, k_n):=\sum^{n-1}_{i=1}k_1\dots k_{i-1}(k_i\odot k_{i+1})k_{i+2}\dots k_n,
\end{equation}
where $k_1,\dots,k_n\in \{x_0,x_1\}$ and $(x_i\odot x_{j}):=\delta_{x_i,x_j}x_i$, for $x_i,x_j\in \{x_0,x_1\}$.
\end{definition}

\begin{definition}[The reduced coaction equation]
	Let $\eta \in \kk$ be a Lie series 
    with $c_{x_0}(\eta)=c_{x_1}(\eta)=0$. The function $r_{\eta}$ associated to $\eta$ is defined to be
\begin{align*}
		r_{\eta}(x)=\sum_{l\ge 0}c_{x^{l+1}_0x_1}(\eta)x^{l+1}
\end{align*}
	and the \textit{reduced coaction equation} is 
	\begin{equation}\label{reduced_coaction_equation}
		\mu(\eta)=-r_{\eta}(x_1)+r_{\eta}(-x_0)-d_0^L(\eta)-d_1^R(\eta).
	\end{equation}
\end{definition}

We denote by $\mathfrak{rc}$ the set of skew-symmetric solutions of the reduced coaction equation, i.e.
\begin{align*}
    \mathfrak{rc}:=\{\psi\in \mathrm{Skew}\mid \text{$\psi$ satisfies \eqref{reduced_coaction_equation}}\}.
\end{align*}

\renewcommand{\thetheorem}{\arabic{theorem}}
\setcounter{theorem}{1}
\begin{lemma}
The vector space $\mathfrak{rc}$ contains the one-dimensional vector space spanned by the element $[x_0,x_1]$.
\end{lemma} 
\begin{proof}
Let $\eta(x_0,x_1) = [x_0,x_1] = x_0x_1 - x_1x_0$. Clearly, $\eta \in \sk$ and using the definitions above, we compute
\begin{align*}
    \mu(\eta) & = \mu(x_0x_1)-\mu(x_1x_0)=0, \quad  r_\eta(x) =x, \\
    d_0^L(\eta) & = -x_1, \quad d_1^R(\eta) = -x_0.
\end{align*}
Thus, 
\begin{align*}
    -r_{\eta}(x_1)+r_{\eta}(-x_0)-d_0^L(\eta)-d_1^R(\eta) = -x_1 + (-x_0) - (-x_1) - (-x_0) = 0 = \mu(\eta),
\end{align*}
so $\eta$ satisfies \eqref{reduced_coaction_equation} and since this condition is linear in $\eta$, the vector space it spans also lies in $\rc$.
\end{proof}

We denote by $\mathfrak{rc}_\lambda$ the subset of $\mathfrak{rc}$ consisting of elements which have coefficient of the commutator $[x_0,x_1]$ equal to $\lambda$;  
we will specifically focus on the case $\lambda=0$. 
We determined the following equivalent polylogarithmic characterizations of $\rc_0$.

\renewcommand{\thetheorem}{\Alph{theorem}}
\setcounter{theorem}{1}
\begin{theorem}[Theorem \ref{th:poly_rc_0}] \label{th:poly_rc_0-C}
If $\psi\in \mathrm{Skew}$ and $c_{x_0x_1}(\psi)=0$, then the following four conditions are equivalent:

    \begin{enumerate}[(i)]
        \item 
        $\psi\in \mathfrak{rc}_0$;
        \item 
        $l^{y,x}_{{\bf{a}},(b_1)}(\alpha)=0, \quad \forall {\bf{a}},(b_1)$;
        \item 
        $l_{{\bf{a}},(b_1)}^{x,y}(\alpha)=0, \quad  \forall {\bf{a}},(b_1)$;
        \item \label{thm:rco-equiv-4}
        $\mu(\psi(-x_0-x_1,x_1))=  d^R_1(\psi(-x_0-x_1,x_1))(x_0+x_1,0) - d ^R_{1}(\psi(-x_0-x_1,x_1))-d^R_1(\psi(-x_0-x_1,x_1))(x_1,0).$  
    \end{enumerate}
\end{theorem}

In Section \ref{sec:linking-dmr-rc}, we prove that $\mathfrak{rc}_0$ is a Lie algebra with respect to the Ihara bracket (see Theorem \ref{theorem:property}). This is done by showing that $\rc_0$ is isomorphic to $\overline{\mathfrak{rc}_0}$, the Lie algebra introduced by \cite{Ren2025} and defined by a variant of the reduced coaction equation \eqref{reduced_coaction_equation} (see \eqref{eq:bar_rc}).
Additionally, we explicitly relate the Lie algebras $\mathfrak{rc}_0$ and $\mathfrak{dmr}_0$; the proof relies on the preceding polylogarithmic descriptions of $\mathfrak{dmr}_0$ and $\mathfrak{rc}_0$. 

\renewcommand{\thetheorem}{\Alph{theorem}}
\setcounter{theorem}{2}
\begin{theorem}[Theorem \ref{theorem:rc_dmr}] \label{theorem:rc_dmr-A}
	Let $\psi\in \mathrm{Skew}$, then the following two conditions are equivalent:
	\begin{enumerate}[(i)]
		\item 
        $\psi\in \mathfrak{dmr}_0$.
		\item 
        $\psi\in \mathfrak{rc}_0$ and for any $\textbf{a,b}$ $\ne (1,\dots,1),(1,\dots,1)$,
			$l^{y,x}_{(a_1,\dots,a_k),(b_1,\dots,b_l)}(\psi_{451}+\psi_{123})=l^{y,x}_{(a_1,\dots,a_k,b_1),(b_2,\dots,b_l)}(\psi_{451}+\psi_{123})$.
	\end{enumerate}
\end{theorem}

Finally, in Section \ref{sec:Kashiwara--Vergne}, we recall the definition of Alekseev-Torossian's \textit{Kashiwara--Vergne Lie algebra} $\krv_2$ and prove the following chain of injective Lie algebra maps, linking the three Lie algebras $\mathfrak{dmr}_0$, $\mathfrak{rc}_0$ and $\mathfrak{krv}_2$: we show that the skew-symmetric elements of $\mathfrak{dmr}_0$ map into $\mathfrak{rc}_0$, and that the elements of $\mathfrak{rc}_0$ satisfying the Vep condition map into $\mathfrak{krv}_2$. 

\renewcommand{\thetheorem}{\Alph{theorem}}
\setcounter{theorem}{3}
\begin{theorem}[Theorem \ref{th:Kashiwara--Vergne}] \label{th:Kashiwara--Vergne-E}
The maps $L$ and $L_1$
\begin{align*}
        L: \mathfrak{dmr}_0 \cap \mathrm{Skew} & \to \mathfrak{rc}_0 \\
        \psi & \mapsto \psi,
    \end{align*}
\begin{align*}
\mathfrak{dmr}_0\cap \mathrm{Skew}\cap\mathrm{Vep}&\xrightarrow{L} \mathfrak{rc}_0\cap \mathrm{Vep}\xrightarrow{L_1} \mathfrak{krv}_2&\\
\psi(x_0,x_1)&\mapsto \psi(x_0,x_1)\mapsto (-\psi(-x_0-x_1, x_0), -\psi(-x_0-x_1, x_1)),
\end{align*}
are injective Lie algebra maps.
\end{theorem}

We close this discussion with a few remarks relating our approach to prior work. Schneps \cite{Schneps2025} and Enriquez-Furusho \cite{EF4} independently proved that $\mathfrak{dmr}_0 \hookrightarrow \krv_2$. As our corresponding map factors through $\rc_0$ and requires the Vep equation, it is open whether it coincides with theirs.
The reduced coaction also appears as $R$ in Kuno's work on emergent braids \cite{Kuno2025}. 
The main difference is the skew-symmetric condition, as highlighted in Remark \ref{rem:skew}; after identifying the notations, our equation \eqref{thm:rco-equiv-4} of Theorem \ref{th:poly_rc_0-C} coincides with equation (21) in \cite{Kuno2025}.  We hope to investigate this relationship more closely in future work. \\

{\bf Acknowledgements.} This paper grew out of a working group with Anton Alekseev, Francis Brown, Florian Naef and Pavol \v{S}evera, held at the University of Geneva from November 2023 to June 2024. We would like to thank them for all the discussions and for introducing us to this interesting topic. We would also like to thank Benjamin Enriquez and Hidekazu Furusho for answering questions about their papers and for their suggestions. We are also grateful to Yusuke Kuno, Dror Bar-Natan, Lucas Dauger and Khalef Yaddaden for sharing their work and insightful discussions. 
The work of M.H. is partially supported by the Swiss NSF grant $200020$-$200400$. The work of M.R. is supported by the SNSF postdoc mobility grant  P$500$PT$\_230340$. Both authors acknowledge the support of NCCR SwissMAP.

\renewcommand{\thetheorem}{\arabic{theorem}}
\setcounter{theorem}{3}

\section{Notations and preliminaries on polylogarithms}
\label{section:notation_polylogarithms}
We begin by reviewing the relevant background on polylogarithms, following the exposition of \cite{Furusho2011}. Let $\textbf{a}=(a_1,\dots,a_k)\in \Z^k_{>0}$, then its \textit{weight} and its \textit{depth} are respectively defined to be $\text{wt}(\textbf{a})=a_1+\dots+a_k$ and $\text{dp}(\textbf{a})=k$. 

In this work, we restrict our attention to multiple polylogarithm functions in \textit{one} and \textit{two} variables,
which take the form:
\begin{equation}
\label{eq:polylogs-1-2}
\begin{split}
	\Li_{(a_1, \ldots, a_k)}(z) & :=\sum_{0<m_1<\dots<m_k}\frac{z^{m_k}}{m^{a_1}_1\dots m^{a_k}_k}; \\
    \Li_{(a_1, \ldots, a_k),(b_1, \ldots, b_l)}(x,y) & :=\sum_{\substack{0<m_1<\dots< m_k\\ <n_1<\dots <n_l}}\frac{x^{m_k}y^{n_l}}{m^{a_1}_1\dots m^{a_k}_kn^{b_1}_1\dots n^{b_l}_l}.
\end{split}
\end{equation}

These functions respectively satisfy the following differential equations:

\begin{equation}
	\frac{d}{dz}\text{Li}_{\bf{a}}(z)=\begin{cases}
		\frac{1}{z}\text{Li}_{(a_1,\dots,a_{k-1},a_k-1)}(z), & \text{if $a_k\ne 1$}\\
		\frac{1}{1-z}\text{Li}_{(a_1,\dots,a_{k-1})}(z), & \text{if $a_k=1,k\ne 1$}\\
		\frac{1}{1-z}, & \text{if } a_k=1, k=1.
	\end{cases}
\end{equation}
and
\begin{equation}\label{equation:2_variables_different_equation}
\begin{split}
	&\frac{d}{dx}\text{Li}_{\bf{a},\bf{b}}(x,y)=\begin{cases}
		\frac{1}{x}\text{Li}_{(a_1,\dots,a_{k-1},a_k-1),{\bf{b}}}(x,y),&~\text{if}~a_k\ne 1 \\
		\frac{1}{1-x}\text{Li}_{(a_1,\dots,a_{k-1}),\bold b}(x,y)-(\frac{1}{x}+\frac{1}{1-x})\text{Li}_{(a_1,\dots,a_{k-1},b_1),(b_2,\dots,b_l)}(x,y)
		&~\text{if}~a_k=1,k\ne 1 ,l\ne 1\\
		\frac{1}{1-x}\text{Li}_{{\bf{b}}}(y)-(\frac{1}{x}+\frac{1}{1-x})\text{Li}_{(b_1),(b_2,\dots,b_l)}(x,y)&~\text{if}~a_k=1,k=1,l\ne 1\\
		\frac{1}{1-x}\text{Li}_{(a_1,\dots,a_{k-1}),{\bf{b}}}(x,y)-(\frac{1}{x}+\frac{1}{1-x})\text{Li}_{(a_1,\dots,a_{k-1},b_1)}(xy)
		&~\text{if}~a_k=1,k\ne 1,l=1\\
		\frac{1}{1-x}\text{Li}_{\bf{b}}(y)-(\frac{1}{x}+\frac{1}{1-x})\text{Li}_{\bf{b}}(xy),&~\text{if}~a_k=1,k=1,l=1
	\end{cases}\\
	&\frac{d}{dy}\text{Li}_{\bf{a},\bf{b}}(x,y)=\begin{cases}
		\frac{1}{y}\text{Li}_{{\bf{a}},(b_1,\dots,b_{l-1},b_l-1)}(x,y) &~\text{if}~b_l\ne 1\\
		\frac{1}{1-y}\text{Li}_{\bf{a},(b_1,\dots,b_{l-1})}(x,y)& ~\text{if}~b_l=1,l\ne 1\\
		\frac{1}{1-y}\text{Li}_{\bf{a}}(xy) &~\text{if}~b_l=1,l=1.
	\end{cases} 
\end{split}
\end{equation}

One important property is the \textit{stuffle relation} (see \cite{Furusho2011, CS}), i.e. 

\begin{align*}
    Li_\textbf{a}(x)Li_\textbf{b}(y)=\sum_{\sigma \in \Sh^{\leq(k,l)}}Li_{\sigma(\textbf{a},\textbf{b})}(\sigma(x,y)),
\end{align*}
where 
\begin{align*}
\Sh^{\leq(k,l)} = \bigcup_{N=1}^{\infty} \Biggl\{ 
       \begin{array}{l|cl}
                      & \sigma \text{ is onto;}\\
             \sigma : \{1, \ldots, k+l\} \to \{1, \ldots, N\} & \sigma(1) < \ldots < \sigma(k);\\
             & \sigma(k+1) < \ldots < \sigma(k+l)
        \end{array}
     \Biggr\},
\end{align*}
and
\begin{equation*}
\begin{split}
    \sigma(x,y) & := \begin{cases}
        xy, \hfill & \text{ if } \sigma^{-1}(N)=k, k+l; \\
        (x,y), \hfill & \text{ if } \sigma^{-1}(N)=k+l; \\
        (y,x), \hfill & \text{ if } \sigma^{-1}(N)=k,
    \end{cases}\\
    \sigma(\textbf{a}, \textbf{b}) & :=((c_1,\ldots, c_j),(c_{j+1}, \ldots, c_N)), 
\end{split}
\end{equation*}
where
    \begin{align*}
        \{j,N\} & =\{\sigma(k), \sigma(k+l)\},\\
        c_i & = \begin{cases}
            a_s + b_{t-k}, \hfill & \text{ if } \sigma^{-1}(i) = \{s,t\} \text{ with } s < t;\\
            a_s, \hfill & \text{ if } \sigma^{-1}(i) = \{s\} \text{ with } s \leq k;\\
            b_{s-k}, \hfill & \text{ if } \sigma^{-1}(i) = \{s\} \text{ with } s > k.\\
        \end{cases}
    \end{align*}

In particular, the functions $\text{Li}_{\bf{a},\bf{b}}(x,y), \text{Li}_{\bf{a},\bf{b}}(y,x),\text{Li}_{\bf{a}}(x),\text{Li}_{\bf{a}}(y)$ and $\text{Li}_{\bf{a}}(xy)$ can be written in the form of Chen's iterated integrals; we use Zhao's form \cite{Zhao} of iterated integrals in two variables.

\begin{ex}\label{ex:dilog}[\cite{Zhao}]  
For example, the double dilogarithm is given by
\begin{equation}\label{eq:diloga}
Li_{(1),(1)}(x,y) = \int_{(0,0)}^{(x,y)} \left( \frac{dx}{1-x} \frac{dy}{1-y} + \left( \frac{dy}{1-y} - \frac{dx}{1-x} - \frac{dx}{x} \right) \frac{d(xy)}{1 - xy} \right).
\end{equation}
\end{ex}

They are functions on the universal covering space of the moduli space of points $\widehat{\mathcal{M}_{0,r+3}}$, where $r$ is the number of variables. We now recall its definition, along with that of the reduced bar construction of these moduli spaces. Let $\mathcal{M}_{0,r+3}$ be the moduli space of $r+3$ different points in $\mathbb{P}^1(k)$ modulo the PGL$_2(k)$ action. It is identified with 
\begin{align*}
\{(z_1,\ldots,z_r)\in (\mathbb{P}^1\setminus \{0,1,\infty\})^r \mid z_i\ne z_j, \text{ for }i\ne j\},
\end{align*}
and under the change of variables
$
x_i:=z_i/z_{i+1} \text{ for } 1\le i\le r-1, x_r:=z_r,
$
we have the identification
\begin{align*}
\mathcal{M}_{0,r+3}\simeq\{(x_1,\ldots,x_r)\in G^{r}_{m} \mid x_i\ne 1, x_{i}x_{i+1}\dots x_{j}\ne 1, 1\le i\le j\le  r\}.
\end{align*}
This coordinate system is called the cubic coordinate system on $\mathcal{M}_{0,r+3}$ and is studied in \cite{Brown2009}. The reduced bar construction
$\mathcal{V}(\mathcal{M}_{0,r+3})$ is the graded dual of the universal enveloping algebra of the pure sphere braid Lie algebra on $r+3$ strands, $\mathfrak{p}_{r+3}$. We now recall its definition.

\begin{definition} \label{def:pure-shere-braid}
    The \textit{pure sphere braid Lie algebra on $r+3$ strands}, $\mathfrak{p}_{r+3}$, is defined by the generators $x_{ij} = x_{ji}$, for $1 \leq i,j \leq r+3$, subject to the relations
    \begin{equation*}
        \begin{split}
            & x_{ii} =0, \quad \forall i \in \{1, \ldots, r+3\}; \\
            & \sum_{j=1}^{r+3} x_{ij} = 0, \quad \forall i \in \{1, \ldots, r+3\}; \\
            & [x_{ij},x_{kl}] = 0, \quad \text{if } \{i,j\} \cap \{k,l\} = \emptyset.
        \end{split}
    \end{equation*}
\end{definition}

\begin{ex}[\cite{Furusho2011}]
We detail the case $r=2$. Let $V_1=H^1_{\rm{DR}}(\mathcal{M}_{0,5})$ and $\widehat{T}(V_1)=\prod_{m\ge 0}(V_1)^{\otimes m}$ be its completed tensor algebra. $V_1$ is a vector space with the basis
\begin{align*}
\omega_{45}:=\frac{dy}{y},\quad \omega_{34}=\frac{dy}{y-1},\quad \omega_{24}:=\frac{d(xy)}{xy-1},\quad \omega_{12}:=\frac{dx}{x},\quad \omega_{23}:=\frac{dx}{x-1}.
\end{align*}

The reduced Bar construction
$\mathcal{V}(\mathcal{M}_{0,5})=\oplus^{\infty}_{m=0} V_m \subset T(V_1)$, where $V_m$ consists of (finite) linear combinations 
\begin{align*}
\sum_{I=(i_m,\ldots,i_1)}c_I\omega_{i_m}\otimes\ldots \otimes \omega_{i_1}\in V_1^{\otimes m}
\end{align*}
which satisfy Chen's integrability condition in $H^1(\mathcal{M}_{0,5})\otimes \ldots H^2(\mathcal{M}_{0,5})\ldots\otimes H^1(\mathcal{M}_{0,5})$, i.e. for all $1\le j<m$,
\begin{align*}
\sum_{I=(i_m,\ldots,i_1)} c_I\omega_{i_m}\ldots \otimes (\omega_{i_{j+1}}\wedge \omega_{i_j})\otimes\ldots \omega_{i_1}=0.
\end{align*}
It is isomorphic as a Hopf algebra to the dual of $U\mathfrak{p}_5$, denoted $(U\mathfrak{p}_5)^{*}$, through the identification of the degree $1$ part $\omega_{45},\omega_{34},\omega_{24},\omega_{12},\omega_{23}$ with $x_{45},x_{34}, x_{24},x_{12}, x_{23}$.
\end{ex}

Let $o$ denote the tangential base point $x_1=x_2=\ldots=x_r=0$ with tangent vector $(1,\ldots,1)$ and let $I_o(\widehat{\mathcal{M}_{0,r+3}})$ denote the homotopy invariant iterated integral on the universal covering space $\widehat{\mathcal{M}_{0,r+3}}$.

\begin{theorem}[\cite{Brown2009}, Equation 3.26]
There is an embedding 
    \begin{align*}
        \rho: \mathcal{V}(\mathcal{M}_{0,r+3}) & \hookrightarrow I_{o}(\widehat{\mathcal{M}_{0,r+3}}) \\
        \sum_{I=(i_m,\ldots,i_1)}c_I[\omega_{i_m} \mid \ldots \mid \omega_{i_1}] & \mapsto \sum_{I}c_I \int_{o}\omega_{i_m}\circ \ldots \circ \omega_{i_1} 
    \end{align*}
and the right-hand side denotes the iterated integral defined by
\begin{align*}
\sum_{I}c_I\int_{0<t_1<\ldots<t_m<1}\omega_{i_m}(\gamma(t_m))\ldots \omega_{i_1}(\gamma(t_1)),
\end{align*}
for all analytic paths $\gamma:(0,1)\to \mathcal{M}_{0,r+3}(\mathbb{C})$ such that $\gamma(0)=(0,\ldots,0)$ and $\dot{\gamma}(0)=(1,\ldots,1).$
\end{theorem}

The elements of the reduced bar construction which correspond to the one- and two-variable multiple polylogarithms defined in \eqref{eq:polylogs-1-2} under the map $\rho$ are denoted by 
\begin{align*}
    l_{\mathbf{a}} :=\rho^{-1}(\Li_{\mathbf{a}}(z)) \in \mathcal{V}(\mathcal{M}_{0,4})
\end{align*}
and 
\begin{align*}
    l^{x,y}_{\bf{a},\bf{b}} :=\rho^{-1}(\Li_{\bold a,\bold b}(x,y)), l^{y,x}_{\bf{a},\bf{b}} :=\rho^{-1}(\Li_{\bold a,\bold b}(y,x)),
    l^{x}_{\bf{a}} :=\rho^{-1}(\Li_{\bold a}(x)), l^{y}_{\bf{a}} :=\rho^{-1}(\Li_{\bold a}(y)),
    l^{xy}_{\bf{a}} :=\rho^{-1}(\Li_{\bold a}(xy)) \in \mathcal{V}(\mathcal{M}_{0,5}).
\end{align*}
These bar words were introduced in \cite{Furusho2011} and are uniquely determined by the differential equations \eqref{equation:2_variables_different_equation}.

\begin{ex}
For $\mathcal{V}(\mathcal{M}_{0,4})$, $\omega_0=\frac{dz}{z}$ and $\omega_1=\frac{dz}{z-1}$, and
\begin{align*}
	l_{\bf{a}}=(-1)^k[\underbrace{\omega_0|\dots|\omega_0|}_{a_k-1}\omega_1|\underbrace{\omega_0|\dots |\omega_0}_{a_{k-1}-1}|\omega_1|\dots |\underbrace{\omega_0|\dots |\omega_0|}_{a_1-1}\omega_1],
\end{align*}
which, evaluated in a series $\varphi=\sum_{W:\text{word}}c_{W}(\varphi)W$, is calculated by 
\begin{align*}
    l_{\bf{a}}(\varphi):=\langle l_{\bf{a}},\varphi \rangle=(-1)^kc_{x^{a_k-1}_0x_1x^{a_{k-1}-1}_0x_1\dots x^{a_1-1}_0x_1}(\varphi).
\end{align*}
\end{ex}

\begin{ex}
The bar word $l^{x,y}_{(1),(1)}\in \mathcal{V}(\mathcal{M}_{0,5})$ corresponding to the double dilogarithms \eqref{eq:diloga} of Example \ref{ex:dilog} is 
\begin{align*}
l^{x,y}_{(1),(1)}= [\frac{dx}{1-x}\,\mid \,\frac{dy}{1-y}]+[\frac{dy}{1-y}-\frac{dx}{x}-\frac{dx}{1-x}\,\mid \,\frac{d(xy)}{1-xy}]
= [\omega_{23}\mid \omega_{34}]-[\omega_{23}-\omega_{34}-\omega_{12}\mid \omega_{24}].
\end{align*}
\end{ex}

\section{Double shuffle relation and defect of the pentagon equation}
\label{sec:poly-descr-double-shuffle}

Let $\psi\in \mathfrak{fr}_k(x_0,x_1)$. We consider the defect $\alpha$ of its pentagon equation
\begin{equation}\label{eq:defect-no-skew}
\alpha\left(\psi(x_0,x_1)\right):=\psi_{451}+\psi_{123}-\psi_{432}-\psi_{215}-\psi_{543},
\end{equation} 
where we use the shorthand $\psi_{ijk} := \psi(x_{ij},x_{jk})\in \mathfrak{p}_5$ for $1\le i,j,k\le 5$.

\begin{rmk}
We consider this specific form of the pentagon equation since we do not assume $\psi$ to be skew symmetric here.
\end{rmk}

In this section, we study the defect of the pentagon equation and its relation to the double shuffle Lie algebra $\mathfrak{dmr}_0$, which was defined in Definition \ref{def:double-shuffle}. The main tool we use here is polylogarithm calculations, which we subdivide into two cases.

\subsection{Case $\textbf{a},\textbf{b} \neq (1,\ldots, 1),(1,\ldots,1)$}
\label{subsection:indexnot1}

We start by collecting some useful results linking polylogarithm functions to the defect of the pentagon equation; they are a compilation of Lemmas $4.1$ and $4.2$ from \cite{Furusho2011} and Lemma $3$ from \cite{CS}.

\begin{lemma}
\label{lemma:polylogs-compilation}
Let $\psi \in \mathfrak{fr}_k(x_0,x_1)$ be a Lie series. Then,
\begin{enumerate}[(i)]
    \item \label{lemma:543}$l^{y,x}_{\bf{a},\bf{b}}(\psi_{543}) =0,\quad \text{for any}~ {\bf{a},\bf{b}};$
    \item \label{lemma:215} $l^{y,x}_{\bf{a},\bf{b}}(\psi_{215}) =l_{\bf{a}\bf{b}}(\psi),\quad \text{for any}~ {\bf{a},\bf{b}}; $
    \item \label{lemma:432} $l^{y,x}_{\bf{a},\bf{b}}(\psi_{432}) =0,\quad \text{for}~ {\bf{a},\bf{b}} \ne (1,\dots,1), (1,\dots,1);$
    \item \label{eq: 5 term_1} $l^{xy}_{\bf{a}}(\psi_{451}+\psi_{123})=l_{\bf{a}}(\psi),\quad \text{for any}~ {\bf{a}};$
    \item \label{eq: 5 term_2} $l^{x,y}_{\bf{a},\bf{b}}(\psi_{451}+\psi_{123})=l_{\bf{a}\bf{b}}(\psi),\quad \text{for any}~ {\bf{a},\bf{b}}.$
\end{enumerate}
\end{lemma}

\begin{lemma}[\cite{Furusho2011}, Equation (3.2)] Let $\varphi\in \mathfrak{p}_5$ be a Lie series, then the following \emph{series shuffle relations (stuffle)} modulo product holds for any index $\bf{a},\bf{b}$,
	\begin{equation}\label{eq:stuffle 2 variable}
		\sum_{\substack{\sigma\in \text{Sh}^{\le(k,l)}, \\
        \sigma^{-1}(N)=k,k+l}} l^{xy}_{\sigma(\bf{a},\bf{b})}(\varphi) 
        +\sum_{\substack{\sigma\in \text{Sh}^{\le(k,l)}, \\ \sigma^{-1}(N)=k+l}} l^{x,y}_{\sigma(\bf{a},\bf{b})}(\varphi) 
        +\sum_{\substack{\sigma\in \text{Sh}^{\le(k,l)}, \\ \sigma^{-1}(N)=k}} l^{y,x}_{\sigma(\bf{a},\bf{b})}(\varphi)=0.
	\end{equation}

\end{lemma}

\begin{proposition} \label{prop:sum}
	If $\psi \in \mathfrak{dmr}_0$, then for any $\bf{a,b}$,
	\begin{equation}\label{eq:sum}
		\sum_{\substack{\sigma\in \text{Sh}^{\le(k,l)}, \\ \sigma^{-1}(N)=k}} \left( l^{y,x}_{\sigma(\bf{a},\bf{b})}(\psi_{451}+\psi_{123})-l_{\sigma(\bf{a},\bf{b})}(\psi)\right)=\sum_{\sigma\in \text{Sh}^{\le(k,l)}}l_{\sigma(\bf{a},\bf{b})}(\psi_{\text{corr}}).
	\end{equation}
\end{proposition}

\begin{pf}
	For $\psi\in \mathfrak{dmr}_0$, by definition of $\Delta_{*}$, the condition $\Delta_{*}(\psi_{*})=1\otimes \psi_{*}+\psi_{*}\otimes 1$ is equivalent to 
	\begin{equation}\label{eq:stuffle 1 variable}
		\sum_{\sigma\in \text{Sh}^{\le(k,l)}}	l_{\sigma(\bf{a},\bf{b})}(\psi_*)=\sum_{\sigma\in \text{Sh}^{\le(k,l)}}	l_{\sigma(\bf{a},\bf{b})}(\psi+\psi_{\text{corr}})=0.
	\end{equation}
	We evaluate equation \eqref{eq:stuffle 2 variable} in the element $\psi_{451}+\psi_{123}$ and compare the result with equation \eqref{eq:stuffle 1 variable}. Finally, we conclude using parts \eqref{eq: 5 term_1} and \eqref{eq: 5 term_2} of Lemma \ref{lemma:polylogs-compilation}. 
%
	
\end{pf}

\begin{proposition}\label{prop:admissible}
	If $\psi \in \mathfrak{dmr}_0$, then for ${\bf{a},\bf{b}}\ne (1,\dots,1),(1,\dots,1)$,
	\begin{equation}\label{eq:2 to 1}
		l^{y,x}_{\bf{a},\bf{b}}(\psi_{451}+\psi_{123})=l_{\bf{ab}}(\psi).
	\end{equation}
\end{proposition}

\begin{pf}
	For the case ${\bf{c},\bf{d}}\ne (1,\dots,1),(1,\dots,1)$, we have that $\sigma(\textbf{c},\textbf{d})\ne (1,\ldots 1)$ for any $\sigma\in \text{Sh}^{\le\left(\dep(\textbf{c}), \dep(\textbf{d})\right)}$, and thus $l_{\sigma(\bf{c},\bf{d})}(\psi_{\text{corr}})=0$ by definition of $\psi_{\text{corr}}$. Applying equation \eqref{eq:sum}, we get
	\begin{equation} \label{eq:proofbis}
		\sum_{\substack{\sigma\in \text{Sh}^{\le(\dep(\textbf{c}), \dep(\textbf{d}))}, \\ \sigma^{-1}(N)=\dep(\textbf{c})}} l^{y,x}_{\sigma(\bf{c},\bf{d})}(\psi_{451}+\psi_{123}) = \sum_{\substack{\sigma\in \text{Sh}^{\le(\dep(\textbf{c}),\dep(\textbf{d}))}, \\ \sigma^{-1}(N)=\dep(\textbf{c})}} l_{\sigma(\bf{c},\bf{d})}(\psi).
	\end{equation}\\
    
    Fix $\dep(\textbf{a})=k$, we proceed by induction on the depth of $\bf{b}$.
	
	\begin{enumerate}
		\item Assume that ${\bf{b}}=(b_1)$ is of depth 1. Take $\bf{c}=\bf{b}$ and $\bf{d}=\bf{a}$ in  \eqref{eq:proofbis}, then 
            $
                |\{\sigma \mid \sigma\in \text{Sh}^{\le(1,k)},\sigma^{-1}(N)=1 \}|=1,
            $
        and this $\sigma$ satisfies $\sigma(\bf{c},\bf{d})=(a,b)$. Thus, equation \eqref{eq:proofbis} implies \eqref{eq:2 to 1} for {$\bf{b}$} of depth $1$.
        
		\item Assume that equation \eqref{eq:2 to 1} is true for any $\bf{b}$ of depth $1,\dots,n$ and suppose that the depth of $\bf{b}$ is $n+1$. We choose $\bf{c}=\bf{b}$ and $\bf{d}=\bf{a}$ in  \eqref{eq:proofbis}, then we have 
		$$
		|\{\sigma \mid \sigma\in \text{Sh}^{\le(n+1,k)},\sigma^{-1}(N)=n+1,{(\bf{e},\bf{f})=\sigma(c,d)}, \dep(\textbf{f})=n+1 \}|=1,
		$$
		and this $\sigma$ satisfies $\sigma(\bf{c},\bf{d})=(\bf{a},\bf{b})$. Since for the shuffle $(\bf{e},\bf{f})=\sigma(c,d)$, $f$ has depth at most $n+1$, we then prove the result in the depth $n+1$ case by using equation \eqref{eq:proofbis} together with the induction hypothesis.
	\end{enumerate}
\end{pf}

\begin{proposition}\label{prop:defect_admissible}
	If $\psi\in \mathfrak{dmr}_0$, then 
	$
		l^{y,x}_{{\bf{a},\bf{b}}}(\alpha)=0\quad \text{for}~ {\bf{a},\bf{b}} \ne (1,\dots,1), (1,\dots,1).
	$
\end{proposition}

\begin{pf}
	Proposition \ref{prop:admissible} and Lemma \ref{lemma:polylogs-compilation}\eqref{lemma:215} imply that $l^{y,x}_{{\bf{a},\bf{b}}}(\psi_{451}+\psi_{123}-\psi_{215})=0$, and the relations in parts \eqref{lemma:543} and \eqref{lemma:432} of Lemma \ref{lemma:polylogs-compilation} finish the proof.
\end{pf}

The following proposition is a Lie algebraic version of \cite[Proposition 3.2]{Furusho2011}; see also \cite[Theorem 3]{CS}. We sketch a proof here using purely polylogarithm calculations.
 
\begin{proposition}\label{prop:condition-dmr0}
Let $\psi\in \mathfrak{fr}_k(x_0,x_1)$ be such that $c_{x_0}(\psi)=c_{x_1}(\psi)= c_{x_0 x_1}(\psi)=0$.  If for ${\bf{a},\bf{b}}\ne (1,\dots,1),(1, \ldots, 1)$
	\begin{equation}\label{non_admissible}
		l^{y,x}_{\bf{a},\bf{b}}(\psi_{451}+\psi_{123})=l_{\bf{a}\bf{b}}(\psi),
	\end{equation}
	then $\psi\in \mathfrak{dmr}_0$.
\end{proposition}

\begin{pf}
	From \eqref{non_admissible}, we know that $l_{\bf{a}}(\psi)$ satisfies the series shuffle relations for all $\bf{a}$ $\neq (1,\dots,1)$. We can extend $\psi$ uniquely to $\psi^S$ which is defined to be 
$$
	l_{\underbrace{(1,\dots,1)}_{n}}(\psi^{S}):=\frac{(-1)^{n}}{n}c_{x^{n-1}_0x_1}(\psi),\quad l_{\bf{a}}(\psi^S):=l_{\bf{a}}(\psi),\quad \text{for } {\bf{a}}\ne (1,\dots,1).
$$
    It is easy to check that $\psi^S$ satisfies the series shuffle relations for all $\bf{a}$, and since
	 $l_{\bf{a}}(\psi_{*})=l_{\bf{a}}(\psi^S)$ for all $\bf{a}$, we conclude $\psi\in \mathfrak{dmr}_0$.
\end{pf}

We conclude the study of the case $\textbf{a},\textbf{b} \neq (1, \ldots, 1),(1, \ldots, 1)$ with the main theorem of this subsection.

\begin{theorem}[Theorem \ref{theorem:dmr_defect-B}]
\label{theorem:dmr_defect}
	Let $\psi\in \mathfrak{fr}_k(x_0,x_1)$ be such that $c_{x_0}(\psi)=c_{x_1}(\psi)= c_{x_0 x_1}(\psi)=0$, then the following two conditions are equivalent:
	\begin{enumerate}[(i)] 
		\item $\psi\in \mathfrak{dmr}_0$;
		\item $l^{y,x}_{\bf{a},\bf{b}}(\alpha)=0,$ for ${\bf{a},\bf{b}}\ne (1,\dots,1),(1,\dots,1)$.
	\end{enumerate}
\end{theorem}
\begin{pf}
	This result follows from Propositions \ref{prop:defect_admissible} and \ref{prop:condition-dmr0}.
\end{pf}

\subsection{Case $\textbf{a},\textbf{b} = (1,\ldots, 1),(1,\ldots,1)$}
\label{subsection:indexis1}
In this subsection, we take the indices to be $\textbf{a},\textbf{b}= \underbrace{(1,\ldots, 1)}_k,\underbrace{(1,\ldots,1)}_l$, i.e. $\dpt(\textbf{a})=k$ and $\dpt(\textbf{b})=l$.

\begin{proposition}\label{prop:nonadmissible1111} 
	If $\psi\in \mathfrak{dmr}_0$, 
    if ${\bf{a}}=\underbrace{(1,\dots,1)}_{k}$, ${\bf{b}}=\underbrace{(1,\dots,1)}_l$, then 
	
	$$
		\sum_{
			\substack{
				\sigma\in \text{Sh}^{\le (k,l)},\sigma^{-1}(N)=k,\\
				\sigma({\bf{a},\bf{b}})=(1,\dots,1),(1,\dots,1)}} l^{y,x}_{\sigma({\bf{a},\bf{b}})}(\psi_{451}+\psi_{123})=(-1)^{k+l}\frac{(l+k-1)!}{k! \hspace{0.5em} l!}c_{x^{k+l-1}_0x_1}(\psi).
	$$
\end{proposition}

\begin{pf}
	For ${\bf{a}}=\underbrace{(1,\dots,1)}_{k}$, ${\bf{b}}=\underbrace{(1,\dots,1)}_l$, we have 
$$
		\sum_{\sigma\in \text{Sh}^{\le(k,l)}}l_{\sigma(\bf{a},\bf{b})}(\psi_{\text{corr}}) = \lvert\{\sigma \mid \sigma({\bf{a},\bf{b}})=\underbrace{(1,\dots,1)}_{k+l}\} \rvert \hspace{0.5em} l_{\underbrace{(1,\dots,1)}_{k+l}}(\psi_{\text{corr}})
		=(-1)^{k+l}\frac{(l+k-1)!}{k! \hspace{0.5em}l!}c_{x^{k+l-1}_0x_1}(\psi).
$$
	Furthermore from Proposition \ref{prop:admissible}, we know that for $\sigma({\bf{a},\bf{b}})\ne (1,\dots,1),(1,\dots,1)$,
	$
		l^{y,x}_{\sigma(\bf{a},\bf{b})}(\psi_{451}+\psi_{123})=l_{\sigma(\bf{a},\bf{b})}(\psi), 
	$
	so the terms in the sum of the LHS of \eqref{eq:sum} cancel out unless $\sigma(\textbf{a}, \textbf{b})=(1, \ldots, 1),(1, \ldots, 1)$. Thus, we get the relation
	$$
		\sum_{
			\substack{
				\sigma\in \text{Sh}^{\le (k,l)},\sigma^{-1}(N)=k,\\
				\sigma({\bf{a},\bf{b}})=(1,\dots,1),(1,\dots,1)}}
                \left(l^{y,x}_{\sigma({\bf{a},\bf{b}})}(\psi_{451}+\psi_{123}) - l_{\sigma(\bf{a},\bf{b})}(\psi)\right) = (-1)^{k+l}\frac{(l+k-1)!}{k! \hspace{0.5em}l!}c_{x^{k+l-1}_0x_1}(\psi).
    $$
    
	Finally, using the fact that $\psi$ is a Lie series  with $c_{x^n_1}(\psi)=0$ for $n\ge 1$, we have
	$$
		\sum_{
			\substack{
				\sigma\in \text{Sh}^{\le (k,l)},\sigma^{-1}(N)=k,\\
				\sigma({\bf{a},\bf{b}})=(1,\dots,1),(1,\dots,1)}
		}l_{\sigma(\bf{a},\bf{b})}(\psi)=0,
	$$
    which proves the formula.
\end{pf}

\begin{rmk}
	Recursively using the above formula, one could determine the formula of $l^{y,x}_{\bf{a},\bf{b}}(\psi_{451}+\psi_{123})$, for indices ${\bf{a,b}}=(1,\dots,1),(1,\dots,1)$ of any depth.
\end{rmk}

\subsubsection{With skew symmetry}
\label{subsubsec:skew}

Before continuing our calculations for the polylogarithms, we first recall a well-known property of the space $\mathrm{Skew}$ defined in \eqref{def:Skew}, for later use.

\begin{lemma}\label{lem:Skew-Lie-subalgebra}
The vector space
\[
\mathrm{Skew}
=
\left\{
\psi\in\mathfrak{fr}^{>1}_k(x_0,x_1)
\;\middle|\;
\psi(x_0,x_1)=-\psi(x_1,x_0)
\right\}
\]
is a Lie subalgebra of
\(\mathfrak{fr}^{>1}_k(x_0,x_1)\) with respect to the Ihara bracket \eqref{eq:Ihara_bracket}.
\end{lemma}

\begin{proof}
We first observe that $\frone$ is stable under the Ihara bracket. 
Let \(\psi_1,\psi_2\in\mathfrak{fr}^{>1}_k(x_0,x_1)\), then by definition, $\psi_j(x_0,0)=\psi_j(0,x_1)=0$ for $j=1,2$ and
the three Lie series $d_{\psi_2}(\psi_1), d_{\psi_1}(\psi_2),
[\psi_1,\psi_2]$ vanish after setting either \(x_0=0\) or \(x_1=0\). Hence, $\{\psi_1,\psi_2\}\in \frone$. Consider $\sigma_{0,1}:A\to A$ the algebra automorphism exchanging \(x_0\) and \(x_1\);
then $\psi\in\mathrm{Skew}$ is equivalent to $\sigma_{0,1}(\psi)=-\psi.$

For \(\eta\in\mathfrak{fr}_k(x_0,x_1)\), define the derivation
\(\overline d_\eta\) by 
\begin{align*}
    \overline d_\eta: x_0 & \mapsto [x_0,\eta] \\
                        x_1 & \mapsto 0.
\end{align*}

A direct calculation on the generators gives $\sigma_{0,1}\circ d_\eta\circ\sigma_{0,1}
= \overline d_{\sigma_{0,1}(\eta)}$. Moreover
$\overline d_\eta=-d_\eta-\operatorname{ad}_\eta$, where
\(\operatorname{ad}_\eta(a)=[\eta,a]\), since both sides send \(x_0\) to \([x_0,\eta]\) and \(x_1\) to \(0\). It follows that
\begin{align*}
\sigma_{0,1}\bigl(\{\psi_1,\psi_2\}\bigr)
={}&
\overline d_{\sigma_{0,1}(\psi_2)}
   \bigl(\sigma_{0,1}(\psi_1)\bigr)
-
\overline d_{\sigma_{0,1}(\psi_1)}
   \bigl(\sigma_{0,1}(\psi_2)\bigr) 
-
\bigl[
\sigma_{0,1}(\psi_1),
\sigma_{0,1}(\psi_2)
\bigr]                                                    \\
={}&
-\left\{
\sigma_{0,1}(\psi_1),
\sigma_{0,1}(\psi_2)
\right\},
\end{align*}
i.e. \(\sigma_{0,1}\) is an anti-automorphism for the Ihara bracket. \\

Let \(\psi_1,\psi_2\in\mathrm{Skew}\), then
$\sigma_{0,1}(\psi_j)=-\psi_j$ for $j=1,2$, and we have
\[
\begin{aligned}
\sigma_{0,1}\bigl(\{\psi_1,\psi_2\}\bigr)
&=-\{-\psi_1,-\psi_2\} =-\{\psi_1,\psi_2\},
\end{aligned}
\]
thus
\(\{\psi_1,\psi_2\}\in\mathrm{Skew}\), which proves the
claim.
\end{proof}

We henceforth assume that $\psi \in \mathrm{Skew}$; in this case, the defect \eqref{eq:defect-no-skew} is equivalent to
$$
\alpha\left(\psi(x_0,x_1)\right)
=\psi_{123}+\psi_{234}+\psi_{345}+\psi_{451}+\psi_{512}.
$$

Let $\sigma$ and $\tau$ be the following generators of the dihedral group $D_5$
\begin{equation}\label{eq:dihedral_generator}
	\tau=(15)(24)(3),\quad \sigma=(12345),
\end{equation}
then it follows from the action of the $D_5$ on $U\mathfrak{p}_5$ and the skew symmetry of $\psi$ that $\alpha$ has dihedral symmetry, that is,
\begin{equation}\label{eq:dihedral_symmetry}
	\alpha^{\sigma}=\alpha,\quad \alpha^{\tau}=-\alpha.
\end{equation}

We introduce the projection maps 
\begin{equation}
\label{eq:strand-projection}
\rm{pr}_i:\mathfrak{p}_5 \to \mathfrak{p}_4,
\end{equation}
which delete a given strand $i=1,\ldots,5$ of $\mathfrak{p}_5$.
They are surjective Lie algebra morphisms, and each $\ker \pr_{i}$ is a free Lie algebra on three generators. Furthermore, the maps 
$
\pr_{i,j} :=(\pr_i,\pr_j) : \mathfrak{p}_5 \mapsto \mathfrak{p}_4\oplus \mathfrak{p}_4,
$
are surjective Lie algebra morphisms, with 
    \begin{equation}
    \label{eq:kerprij}
        \ker \pr_{i,j} =\ker \pr_{i}\cap \ker \pr_j=(x_{ij}),
    \end{equation}
where $(x_{ij})$ is a Lie ideal in $\mathfrak{p}_5$. This last fact follows from $\ker \pr_i$ and $\ker \pr_j$ being free Lie algebras.

\begin{ex}[\cite{EF1}, Section 5.1.2]
Recall that $\mathfrak{p}_4\simeq \mathfrak{fr}_k(x_0,x_1)$, where $x_0=~x_{14}=~x_{23}$, $x_{1}=x_{12}=x_{34}$ and $x_{\infty}:=-x_0-x_1=x_{13}=x_{24}$. The map $\rm{pr}_2 : \mathfrak{p}_5\to \mathfrak{p}_4$ is explicitly given by
\begin{align*}
    x_{12}, x_{23}, x_{24}, x_{25} \mapsto 0; \quad 
    x_{14}, x_{35} \mapsto x_0; \quad 
    x_{15}, x_{34} \mapsto x_1; \quad
    x_{13}, x_{45} \mapsto x_\infty.
\end{align*}

Moreover, $\ker \rm{pr}_2$ is freely generated by $x_{12},x_{23},x_{24}$, which coincides with the Lie subalgebra of $\mathfrak{p}_5$ generated by $x_{12},x_{23},x_{24},x_{25}$. Note that here, we have relabelled $(1,3,4,5)$ as $(1,3,4,2)$. For the remaining projections, the corresponding results follow by permuting the indices.
\end{ex}

\begin{lemma}\label{lemma:skew}
If $\psi \in$ Skew, then 
\begin{enumerate}[(i)]
\item \label{lemma:skew-i}
    $\alpha\in \ker \pr_1\cap \ker \pr_2\cap\ker \pr_3\cap\ker  \pr_4\cap\ker \pr_5$; 
\item \label{lemma:skew-ii}
    $\alpha\in (x_{ij})$, for $i\ne j \in \{1,2,3,4,5\}$.
\end{enumerate}
\end{lemma}
\begin{proof}
As $\alpha$ is cyclically symmetric, it suffices to check the result for one index. For example,
$$
\pr_2(\alpha)=\psi_{451}+\psi_{543}=\psi(x_{13},x_{34})-\psi(x_{13},x_{34})=0.
$$

It follows that $\rm{pr}_{i,j}(\alpha)=0$, and equivalently $\alpha\in \ker\rm{pr}_{i,j}=(x_{ij}).$
\end{proof}

\begin{lemma}
\label{lemma:l-kerpr2}
For all indices ${\bf{a,b}}=(a_1, \ldots, a_k),(b_1, \ldots, b_l)$,
\begin{equation}
    \label{eq:l-kerpr2}
	\eval{l^{y,x}_{{\bf{a}},{\bf{b}}}}_{\ker \pr_2}=(-1)^{k+l}\omega^{b_l-1}_{12}\omega_{23}\dots \omega_{12}^{b_1-1}\omega_{23}\omega^{a_k-1}_{12}\omega_{24}\dots \omega^{a_1-1}_{12}\omega_{24},
\end{equation}
i.e. for any $\alpha\in \ker \rm{pr}_2$, 
$
l^{y,x}_{\bold a,\bold b}(\alpha)=(-1)^{k+l}\omega^{b_l-1}_{12}\omega_{23}\dots \omega_{12}^{b_1-1}\omega_{23}\omega^{a_k-1}_{12}\omega_{24}\dots \omega^{a_1-1}_{12}\omega_{24}(\alpha).
$
\end{lemma}

\begin{proof}
	We look at the words of $l^{y,x}_{{\bf{a}},{\bf{b}}}$ that contain only $\omega_{12},\omega_{24}$ and $\omega_{23}$. It suffices to compute the restriction of $l^{x,y}_{\bf{a},\bf{b}}$ to $\ker \pr_4$; it satisfies the differential equation-induced recursive relations:
	\begin{itemize}
		\item if $b_l>1,l\ge1$, $l^{x,y}_{{\bf{a}},{\bf{b}}}=(-1)[\omega_{45} \mid l^{x,y}_{{\bf{a}},(b_1,\dots,b_l-1)}]$;
		\item if $b_l=1,l>1$, $l^{x,y}_{{\bf{a}},{\bf{b}}}=(-1)[\omega_{34} \mid l^{x,y}_{{\bf{a}},(b_1,\dots,b_{l-1})}]$;
		\item if $b_l=1,l=1$, $l^{x,y}_{{\bf{a}},{\bf{b}}}=(-1)[\omega_{34} \mid l^{xy}_{{\bf{a}}}]$.
	\end{itemize}
	From those differential equations, we deduce that 
	\begin{align}
    \label{eq:explicit-lxyab}
		l^{x,y}_{\bf{a},\bf{b}} & =(-1)^{k+l}[\omega^{b_l-1}_{45}\mid\omega_{34}\mid\dots \mid\omega^{b_1-1}_{45}\mid\omega_{34}\mid\omega^{a_k-1}_{45}\mid\omega_{24}\mid\dots \mid\omega^{a_1-1}_{45}\mid\omega_{24}], \\
        \label{eq:explicit-lxyab-2}
		l^{y,x}_{\bf{a},\bf{b}} & =(-1)^{k+l}[\omega^{b_l-1}_{12}\mid\omega_{23}\mid\dots \mid\omega^{b_1-1}_{12}\mid\omega_{23}\mid\omega^{a_k-1}_{12}\mid\omega_{24}\mid\dots \mid\omega^{a_1-1}_{12}\mid\omega_{24}].
	\end{align}

It follows that $l^{y,x}_{\bf{a},\bf{b}}$ contains exactly depth $\bf{b}$ $\omega_{23}$'s.
	
\end{proof}

\begin{proposition}\label{prop:vanishing}
	If $\psi \in$ Skew, then
	$
		l^{y,x}_{(1,\dots,1),(1,\dots,1)}(\alpha)=0.
	$
\end{proposition}

\begin{proof}
By Lemma \ref{lemma:skew}, we know
$
	l^{y,x}_{{\bf{a}},{\bf{b}}}(\alpha)= \eval{l^{y,x}_{\bf{a},\bf{b}}}_{{\ker \rm{pr}_2}}(\alpha).
$
When restricted to the case $\textbf{a,b} = (1,\dots,1),(1,\dots,1)$, the expression \eqref{eq:l-kerpr2} of Lemma \ref{lemma:l-kerpr2} becomes
$
	l^{y,x}_{(1,\dots,1),(1,\dots,1)}=(-1)^{k+l}\omega^{l}_{23}\omega^k_{24},
$
which contains no $\omega_{12}$. The result follows because $\alpha\in (x_{12})$, by Lemma \ref{lemma:skew}.
\end{proof}

Together with Theorem \ref{theorem:dmr_defect}, which handled the case when \textbf{a,b} $\neq (1,\ldots,1),(1,\ldots,1)$, we have the following result.
\begin{theorem}\label{th:defect_skew}
	If $\psi \in \mathfrak{dmr}_0\cap \mathrm{Skew}$, then 
	$
		l^{y,x}_{{\bf{a},\bf{b}}}(\alpha)=0,\quad \forall {\bf{a},\bf{b}}.
	$
\end{theorem}

\begin{proof}
It follows from Theorem \ref{theorem:dmr_defect} that 
$
l^{y,x}_{\bf{a},\bf{b}}(\alpha)=0$ for ${\bf{a},\bf{b}}\ne (1,\dots,1),(1,\dots,1),
$
and from Proposition \ref{prop:vanishing} that 
$l^{y,x}_{\bf{a},\bf{b}}(\alpha)=0$ for ${\bf{a},\bf{b}}= (1,\dots,1),(1,\dots,1).$
\end{proof}

As a corollary, we get a similar result to Proposition \ref{prop:nonadmissible1111}, but for the term $\psi_{432}$ of the defect.

\begin{corollary}
Let $\psi\in \mathfrak{dmr}_0\cap \mathrm{Skew}$. If ${\bf{a}}=\underbrace{(1,\dots,1)}_{k}$, ${\bf{b}}=\underbrace{(1,\dots,1)}_l$, then
	$$
		\sum_{
			\substack{
				\sigma\in \text{Sh}^{\le (k,l)},\sigma^{-1}(N)=k,\\
				\sigma({\bf{a},\bf{b}})=(1,\dots,1),(1,\dots,1)}} l^{y,x}_{\sigma({\bf{a},\bf{b}})}(\psi_{432})=(-1)^{k+l}\frac{(l+k-1)!}{k! \hspace{0.5em} l!}c_{x^{k+l-1}_0x_1}(\psi).
	$$
\end{corollary}

\begin{proof}
By Theorem \ref{th:defect_skew}, we know that
$
l^{y,x}_{\bf{a},\bf{b}}(\psi_{451}+\psi_{123}-\psi_{432}-\psi_{215}-\psi_{543})=0,
$
which implies by Lemma \ref{lemma:polylogs-compilation} that $l^{y,x}_{\bf{a},\bf{b}}(\psi_{451}+\psi_{123}-\psi_{432})=l_{\bf{a}\bf{b}}(\psi)$. As $\psi$ is a Lie series, $l_{(1, \ldots, 1)}(\psi) = 0$, and the result follows from Proposition \ref{prop:nonadmissible1111}.
\end{proof}

\begin{rmk}
	Recursively using the above formula, one could determine the formula for $l^{y,x}_{\bf{a},\bf{b}}(\psi_{432})$, for indices ${\bf{a,b}}~=~(1,\dots,1)~,(1,\dots,1)$ of any depth.
\end{rmk}

\section{Polylogarithmic description of the reduced coaction equation}
\label{sec:poly-descr-reduced-coaction}


The aim of this section is to prove the following equivalent characterizations of the reduced coaction Lie algebra $\rc_0$. Note that we consider the defect $\alpha$ to have the form 
    $\psi_{451} + \psi_{123} - \psi_{432} - \psi_{215} - \psi_{543}$,
as in \eqref{eq:defect-no-skew}.

\begin{theorem}[Theorem \ref{th:poly_rc_0-C}]
\label{th:poly_rc_0}
If $\psi\in \mathrm{Skew}$ and $c_{x_0x_1}(\psi)=0$, then the following 4 conditions are equivalent:

    \begin{enumerate}[(i)]
        \item \label{thm:rco-equiv-1}
        $\psi\in \mathfrak{rc}_0$;
        \item \label{thm:rco-equiv-2}
        $l^{y,x}_{{\bf{a}},(b_1)}(\alpha)=0, \quad \forall {\bf{a}},(b_1)$;
        \item \label{thm:rco-equiv-3}
        $l_{{\bf{a}},(b_1)}^{x,y}(\alpha)=0, \quad  \forall {\bf{a}},(b_1)$;
        \item \label{thm:rco-equiv-4}
        $\mu(\psi(-x_0-x_1,x_1))=  d^R_1(\psi(-x_0-x_1,x_1))(x_0+x_1,0) - d ^R_{1}(\psi(-x_0-x_1,x_1))-d^R_1(\psi(-x_0-x_1,x_1))(x_1,0).$  
    \end{enumerate}
\end{theorem}

We prove the equivalence of \eqref{thm:rco-equiv-1} and \eqref{thm:rco-equiv-2} in Subsection \ref{subsection:3_1}, the equivalence of \eqref{thm:rco-equiv-3} and \eqref{thm:rco-equiv-4} in \ref{subsection:3_2} and finally, the equivalence of \eqref{thm:rco-equiv-2} and \eqref{thm:rco-equiv-3} in Subsection \ref{subsection:3_3}.

\subsection{Fox pairing and $2$-cocycle} \label{subsection:3_1}
Recall that $A$ denotes the Hopf algebra of formal non-commutative power series in $x_0$ and $x_1$, as in \eqref{eq:def-hopf-algebra-A}.
In what follows, we denote by $M$ the space $M \cong A$, endowed with the $A \otimes A$-bimodule structure given by
	\begin{equation}\label{eq:bimodule}
	   (f \otimes g) a (h \otimes k) = \varepsilon(f)\varepsilon(k) \, gah.
	\end{equation}
This bimodule structure plays a key role in this section. The notions of \textit{Fox derivative} and \textit{Fox pairing} were introduced by Massuyeau-Turaev in \cite{Massuyeau2014} as a noncommutative version of the loop operations on a surface. We begin by recalling the tools developed in \cite{ANR} for studying the group version of the reduced coaction equation.

\begin{definition}[Fox derivative]
	A $k$-linear map $\partial: A\to A$ is called a \textit{left Fox derivative} if
	$\partial(ab)=a\partial(b)+\partial(a)\varepsilon(b),$
	for all $a,b\in A$, and a \textit{right Fox derivative} if
	$\partial(ab)=\partial(a)b+\varepsilon(a)\partial(b),
	$ for all $a,b\in A$.
\end{definition}

\begin{ex}[\cite{ANR}, Example $3.2$]
\label{example:left_right fox deri}
	Recall that for any $x\in A$, there are unique presentations
	$$
    x =\varepsilon(x)+ x_0 d^R_0(x) + x_1 d^R_1(x)
    =\varepsilon(x)+ d^{L}_0(x) x_0 + d^{L}_1 (x)x_1.
    $$
	The operation $d^{R}_i$ is a right Fox derivative and $d^L_i$ is a left Fox derivative.
\end{ex}

\begin{definition}[Fox pairing]
    A bilinear map $\rho: A\otimes A\to A$ is a \textit{Fox pairing} if it is a left Fox derivative in its first argument and a right Fox derivative in its second argument. That is,  for all $a,b,c \in A$,
	\begin{equation}\label{eq:product_rule}
		\rho(ab, c) =a\rho(b, c)+\rho(a, c)\varepsilon(b) \quad \text{ and } \quad
        \rho(a,bc) =\rho(a,b)c+\varepsilon(b)\rho(a,c).
	\end{equation}
\end{definition}

\begin{ex}[\cite{ANR}, Example $3.4$]
 We define
    $
        \rho_{\rm{KKS}}~(x_i,x_j)~:=~\delta_{ij}x_i
    $
on the free algebra $\kk$. It uniquely extends to a Fox pairing on $A$ by $\rho_{\rm{KKS}}(x,1)=\rho_{\rm{KKS}}(1,x)=0$, for all $x \in A$ and the product rule \eqref{eq:product_rule}.
\end{ex}
\begin{rmk}
This Fox pairing induces the Kostant-Kirillov-Souriau (KKS) Poisson bracket on the representation space.
\end{rmk}

\begin{proposition}[\cite{ANR}, Proposition 3.5]\label{prop:forpairing2cocycle}
	Let $\rho \colon A \otimes A \to A$ be a Fox pairing. Then,
	$$
	c(a_1 \otimes b_1, a_2 \otimes b_2) := \varepsilon(a_1) \rho(b_1,a_2) \varepsilon(b_2)
	$$
	defines a $2$-cocycle $(A \otimes A)^{\otimes 2} \to M$.
	Moreover, the formula
	$$
		(a_1 \otimes b_1 \oplus c_1) \cdot_{\rho} (a_2 \otimes b_2 \oplus c_2) := a_1 a_2 \otimes b_1 b_2 \oplus \big( \varepsilon(a_1) \, b_1 c_2 + c_1 a_2 \, \varepsilon(b_2) + \varepsilon(a_1) \rho(b_1,a_2) \varepsilon(b_2) \big)
		\label{eqn:multiplicationfromfoxp}
	$$
	defines an algebra structure on $A \otimes A \oplus M$.
\end{proposition}	

\begin{proof}
The associativity of the algebra map is checked explicitly, and it coincides with the $2$-cocycle property.
\end{proof}

Recall that the \textit{Drinfeld-Kohno Lie algebra} $\mathfrak{t}_n$ is the Lie algebra defined by the generators $t_{ij}$ for $1\le i,j\le n$ and the relations 
$$
t_{ij}=t_{ji}; \quad [t_{ij},t_{kl}]=0 \text{ if } \lvert \{i,j,k,l\} \rvert =4 ; \quad [t_{ij}+t_{ik},t_{jk}]=0 \text{ if } \lvert \{i,j,k\} \rvert =3.
$$
The Lie algebra $\mathfrak{t}_n$ has a one-dimensional center $Z(\mathfrak{t}_n)$ and it holds $\mathfrak{t}_n/Z(\mathfrak{t}_n)\simeq \mathfrak{p}_{n+1}$. The map $K_4 : \mathfrak{t}_4 \to \mathfrak{p}_5$ identifies $ \ker \pr_i\subset \mathfrak{t}_4$ with $ \ker \pr_i\subset \mathfrak{p}_5$ for $i \in \{1,2,3,4\}$. In what follows, we only consider the defect $\alpha(\psi)$ which falls inside $\ker \pr_i$; for convenience, we will use the same notation $x_{ij}$ for the corresponding generators $t_{ij}$, for $1\le i<j\le n$.  \\

Define $U\mathfrak{p}^{i,j}_5 := \ker (\pr_i\circ \pr_j)$, for $i\ne j$, where $\pr_i$ and $\pr_j$ are projections in $U\mathfrak{t}_4$ corresponding to those of \eqref{eq:strand-projection}. Note that the superscript $i,j$ is taken with respect to $U\mathfrak{t}_4$ (rather than $U\mathfrak{t}_3$). For example,
$$
\pr_2\circ\pr_3: U\mathfrak{t}_4\xrightarrow{\pr_j} U\mathfrak{t}_3\xrightarrow{\pr_i}k[x_{14}].
$$
Moreover, denote by $I^{ij}:=(x_{ij})$ the ideal generated by $x_{ij}$ in $U\mathfrak{p}^{i,j}_5$.


\begin{lemma}[\cite{ANR}, Proposition $5.1$]\label{eq:4T_relations}
$U\mathfrak{p}^{2,3}_5$ admits the presentation
	\begin{equation*}
	U\mathfrak{p}^{2,3}_5=\langle x_{12},x_{23},x_{24},x_{34},x_{13} \mid R^{2,3} \rangle,
    \end{equation*}
where $R^{2,3}$ corresponds to the relations
\begin{equation*}
\label{eq:4T-relations}
\begin{split}
	[x_{13},x_{23}] &=[x_{23},x_{12}], \quad [x_{34},x_{23}]=[x_{23},x_{24}], \quad [x_{13},x_{24}]=0, \\
	[x_{12},x_{13}] &=-[x_{12},x_{23}], \quad [x_{34},x_{24}]=[x_{24},x_{23}], \quad [x_{12},x_{34}]=0.
\end{split}
\end{equation*}	
\end{lemma}

\begin{proposition}[\cite{ANR}, Proposition $5.2$] \label{prop:pi-2-3}
	The assignment
	\begin{align*}
		\pi^{2,3} \colon U\mathfrak{p}^{2,3}_5 &\to (A \otimes A \oplus A, \cdot_{\rho_{\rm_{KKS}}}) \\
		x_{12} &\mapsto x_0 \otimes 1, \quad x_{24}\mapsto x_1\otimes 1,  \\
		x_{13} &\mapsto 1 \otimes x_0,\quad x_{34}\mapsto 1\otimes x_1,\hskip 0.3cm x_{23} \mapsto -e,
	\end{align*}
	where $e := 0 \oplus 1 \in A \otimes A \oplus M$, extends to an algebra homomorphism. Moreover, the map $\pi^{2,3}$ factors through the canonical projection $U\mathfrak{p}^{2,3}_{5} \to U\mathfrak{p}^{2,3}_5/(I^{2,3})^2$.
\end{proposition}

\begin{proof}
The map is defined explicitly on the generators, so it is easy to show that $\pi^{2,3}(R^{2,3})=0$. Moreover, $\pi^{2,3}((I^{2,3})^2)\subset M\cdot_{\rho_{\rm{KKS}}}M=0$.
\end{proof}

We define the following maps, obtained by composing $\pi^{2,3}$ with the natural projections 
$A \otimes A \oplus M \to A \otimes A$ and $A \otimes A \oplus M \to  M$, 
$$
\pi^{2,3}_0 \colon U\mathfrak{p}^{2,3}_5 \to A \otimes A, \hskip 0.3cm
\pi^{2,3}_1 \colon U\mathfrak{p}^{2,3}_5 \to M.
$$	

Since $A\otimes A$ is identified with $k\langle\langle x_{12},x_{24}\rangle\rangle\otimes k\langle\langle x_{13},x_{34}\rangle\rangle$ via $\pi_0^{2,3}$, $I^{2,3}/(I^{2,3})^2$ is an $A\otimes A$ bimodule generated by $x_{23}$. This bimodule structure is given by two-sided multiplication.

\begin{proposition}
\label{prop:projection}
The map $\pi^{2,3}_0$ descends to an isomorphism 
$
        \pi^{2,3}_0:U\mathfrak{p}^{2,3}_5/I^{2,3} \overset{\simeq}{\longrightarrow} A\otimes A.
$
Furthermore, it holds on $I^{2,3}$ that $\pi^{2,3}_1=F^{2,3}\circ (\pi_1^p)^{2,3}$, where
    \begin{equation*}
        \begin{split}
            (\pi_1^p)^{2,3} & :I^{2,3}\to I^{2,3}/(I^{2,3})^2 \text{ is the projection map}, \\
            F^{2,3} & :I^{2,3}/(I^{2,3})^2 \to M \text{ is the bimodule map determined by } F^{2,3}(x_{23})=-1.
        \end{split}
    \end{equation*}
\end{proposition}

The following split sequence of $A\otimes A$ bimodules follows from Proposition \ref{prop:projection}
\begin{equation}\label{eq:exact sequence}
\begin{tikzcd}
    0\arrow{r} & I^{2,3}/(I^{2,3})^2\arrow{r} &U\mathfrak{p}^{2,3}_5/(I^{2,3})^2\arrow{r} &  U\mathfrak{p}^{2,3}_5/I^{2,3}\arrow{r}\arrow[bend left=33]{l}{(\pi^{2,3}_0)^{-1}} & 0,
\end{tikzcd}
\end{equation}
and using the section map $(\pi^{2,3}_{0})^{-1}$, we can write 
\begin{equation*}
\begin{split}
	U\mathfrak{p}^{2,3}_5/(I^{2,3})^2  & \simeq A\otimes A\oplus I^{2,3}/(I^{2,3})^2 \\
    a  & \mapsto (\pi^{2,3}_0(a),(\pi_1^p)^{2,3}(a)).
\end{split}
\end{equation*}

It is convenient to note that if $\psi \in$ 
Skew, the defect $\alpha$ (as in \eqref{eq:defect-no-skew}) takes the equivalent forms
\begin{align}
\alpha(\psi) & =\psi_{123}+\psi_{234}+\psi_{345}+\psi_{451}+\psi_{512} \nonumber \\
& \label{eq:defect-skew-2}
=\psi_{123}+\psi_{234}-\psi(x_{13}+x_{23},x_{34}) +\psi(x_{12}+x_{13},x_{24}+x_{34})-\psi(x_{12},x_{23}+x_{24}),
\end{align}
where the last equality uses the facts that $\psi$ is a Lie series without degree 1 terms, and the relations 
    $x_{51}~=~x_{23}~+~x_{24}~+x_{34}$ and $x_{45} = x_{12} + x_{13} + x_{23}$ of $\mathfrak{p}_5$.
Moreover, the defect is a linear combination of the following coface maps from $k\langle\langle x_{12},x_{23}\rangle\rangle$ to $U\mathfrak{p}_5$:
\begin{equation}
\label{eq:coface-maps-1}
\begin{split}
c_{1,2,3}:\quad&\psi(x_{12},x_{23})\mapsto \psi(x_{12},x_{23})\\
c_{2,3,4}:\quad&\psi(x_{12},x_{23})\mapsto \psi(x_{23},x_{34})\\
c_{12,3,4}:\quad&\psi(x_{12},x_{23})\mapsto \psi(x_{13}+x_{23},x_{34})\\
c_{1,23,4}:\quad&\psi(x_{12},x_{23})\mapsto \psi(x_{12}+x_{13},x_{24}+x_{34})\\
c_{1,2,34}:\quad&\psi(x_{12},x_{23})\mapsto \psi(x_{12},x_{23}+x_{24}).
\end{split}
\end{equation}

\begin{lemma}
\label{lemma:pi_0}
Let $\psi \in \mathfrak{fr}_k(x_0,x_1)$ and
$\alpha(\psi)=\psi_{123}+\psi_{234}+\psi_{345}+\psi_{451}+\psi_{512}.$
If $\psi \in$ Skew, then $\pi^{2,3}_0(\alpha)=0.$
\end{lemma}

\begin{proof}
Suppose that $\psi \in$ Skew, then $\alpha\in (x_{23})$ by Lemma \ref{lemma:skew} and therefore $\pi^{2,3}_0(\alpha)=0$. 

\end{proof}


\begin{lemma}[\cite{ANR}, Propositions 5.8, 5.10]
\label{lem:cab_23} If $\psi\in \mathrm{Skew}$, then
\begin{align*}
& \pi^{2,3}_1\circ c_{1,2,34}(\psi)=-d^R_1(\psi);\\
& \pi^{2,3}_1\circ c_{12,3,4}(\psi)=-d^L_0(\psi); \\
& \pi^{2,3}_1\circ c_{1,23,4}(\psi)=\mu(\psi);\\
& \pi^{2,3}_1\circ c_{2,3,4}(\psi)=r_{\psi}(x_1);\\
& \pi_1^{2,3}\circ c_{1,2,3}(\psi)=-r_{\psi}(-x_0).
\end{align*}
\end{lemma}

\begin{proof}
We start by showing the details for $\pi^{2,3}_1\circ c_{2,3,4}$. Notice that by Proposition \ref{prop:projection}, $\pi_1^{2,3}=F^{2,3}\circ (\pi^p_1)^{2,3}$ and it factors through $I^{2,3}/(I^{2,3})^2$, so it suffices to look at the terms which only contain one $x_{23}$, $\pi^{2,3}_0(x_{34})=1\otimes x_1$ acts trivially from the right, $x_{34}$ can only appear to the left of $x_{23}$. Hence, we have
$$
\pi^{2,3}_1(\psi(x_{23},x_{34}))=-d^L_0(\psi)(0,x_1)=-\sum_{n\ge 0}c_{x^n_1x_0}(\psi)x^n_1,
$$
which, by skew symmetry of $\psi$, is equal to $\sum_{n\ge 0}c_{x^n_0x_1}(\psi)x^n_1=r_{\psi}(x_1)$. Similarly, for $\pi_1^{2,3}\circ c_{1,2,3}$, we have
$$
\pi^{2,3}_1(\psi(x_{12},x_{23}))=-d^R_1(\psi)(x_0,0)=-\sum_{n\ge 0}c_{x_1x^n_0}(\psi)x^n_0
$$
and because $\psi$ is a Lie series, $S(\psi)=-\psi$, so we get $-\sum_{n\ge 0}c_{x_1x^n_0}(\psi)x_0^n=-r_{\psi}(-x_0)$. Note that the terms at $n=0$ vanish since $\psi \in \fr_k^{>1}$.

The other three equations are proved similarly, using the bimodule structure $\pi_1^{2,3}$ and comparing them to the product rule of the Fox pairing and Fox derivative \eqref{eq:product_rule}; for more details, see \cite{ANR}.
\end{proof}

Let $\gamma$ be a Lie series in the ideal $(x_{23})=\ker \rm{pr}_2\cap \ker \rm{pr}_3\subset \ker \rm{pr}_2$. As the number of $x_{23}$'s is a grading of the vector space $\ker \pr_2=\mathfrak{fr}_k(x_{12},x_{24},x_{23})$, we can write $\gamma=\sum_{i\in \N}\gamma^{(i)}$, where $\gamma^{(i)}$ is the $i$-th $x_{23}$-graded component. As $\gamma$ is a Lie series, we have the following lemma.

\begin{lemma}
\label{lemma:gamma-1-P}
Let $\gamma$ be as above, then
$
\gamma^{(1)}=P(\ad_{x_{12}},\ad_{x_{24}})x_{23}
$
for some formal polynomial $P\in k\langle \langle x_0,x_1\rangle\rangle.$
\end{lemma}

We denote by $V_{\rm{kv}}$ the vector subspace of $I^{2,3}$ spanned by the noncommutative polynomials $P(\ad_{x_{12}},\ad_{x_{24}})x_{23}$, where $P(x_0,x_1) \in k\langle \langle x_0,x_1 \rangle \rangle$, and
by $V_{\kv}^Y$ the vector subspace spanned by those $P(x_0,x_1)$ in which every monomial contains at least one occurrence of the second variable $x_1$.

\begin{lemma}\label{lemma:F_injective}
$\eval{F^{2,3}}_{V_{\rm{kv}}}$ is injective.
\end{lemma}
\begin{proof}
$F^{2,3}$ maps $P(\ad_{x_{12}},\ad_{x_{24}})x_{23}$ to $-S(P(x_0,x_1))$. 
\end{proof}

\begin{proposition}\label{prop:kernel and rc}
If $\psi\in \mathrm{Skew}$, then the following are equivalent:
\begin{enumerate}[(i)] 
    \item \label{prop:kernel and rc-1}
    $\alpha(\psi)\in \ker((\pi^{p}_1)^{2,3})$;
    \item  \label{prop:kernel and rc-2}
    $\psi\in \mathfrak{rc}_0$.
\end{enumerate}
\end{proposition}

\begin{proof}
First, observe that $\alpha\in \ker (\rm{pr}_2)$ and its first order graded component with respect to $x_{23}$ is 
\begin{equation}\label{eq:alpha-graded-component-1}
    \alpha^{(1)}=(\pi^p_1)^{2,3}(\alpha)\in V^{Y}_{\rm{kv}};
\end{equation}
this is because $\alpha\in (x_{24})$ by Lemma \ref{lemma:skew}. 
From \eqref{prop:kernel and rc-1} to \eqref{prop:kernel and rc-2}, suppose that $(\pi^p_1)^{2,3}(\alpha)=0$. It follows that 
$
        \pi^{2,3}_1(\alpha)=F^{2,3}\circ (\pi^p_1)^{2,3}(\alpha)=0,
$
that is, using the form \eqref{eq:defect-skew-2} of the defect:
$$
\pi^{2,3}_1\big(\psi_{123}+\psi_{234}+\psi(x_{12}+x_{13},x_{24}+x_{34})-\psi(x_{13}+x_{23},x_{34})-\psi(x_{12},x_{23}+x_{24})\big)=0.
$$
By Lemma \ref{lem:cab_23},
\begin{align*}
   \pi_1^{2,3}(\alpha) & =-r_{\psi}(-x_0)+r_{\psi}(x_1)+\mu(\psi)+d^R_1(\psi)+d^L_0(\psi) = 0\\
   & \iff \mu(\psi) = r_{\psi}(-x_0) - r_{\psi}(x_1) - d^R_1(\psi) - d^L_0(\psi),
\end{align*}
which is exactly \eqref{reduced_coaction_equation}, hence $\psi\in \mathfrak{rc}_0$.

From \eqref{prop:kernel and rc-2} to \eqref{prop:kernel and rc-1}, suppose that $\psi \in \mathfrak{rc}_0$, then we know by  Lemma \ref{lem:cab_23} that $\pi^{2,3}_1(\alpha)=0$. Since $\alpha^{(1)} \in V^Y_{\kv} \subseteq V_{\kv}$ and $\eval{F^{2,3}}_{V_{\rm{kv}}}$ is injective by Lemma \ref{lemma:F_injective}, it follows that $\alpha\in \ker (\pi^p_1)^{2,3}$.
\end{proof}

We now study the restriction of the polylogarithms $l^{y,x}_{\bf{a},\bf{b}}$ to the subspace $\ker \rm{pr}_2$. Note that in this case, the depth of $\bf{b}$ is a grading by Lemma \ref{lemma:l-kerpr2}.

\begin{lemma}\label{lemma:filtration}
Let $\textbf{a},\textbf{b}=(a_1, \ldots, a_k),(b_1, \ldots, b_l)$, then 
	$
		\eval{l^{y,x}_{\bf{a},\bf{b}}}_{\ker \pr_2}\left( (I^{2,3})^t\cap \ker \pr_2 \right) =0, \text{ for } t> \dpt(\bf{b}).
	$
\end{lemma}

\begin{proof}
By Lemma \ref{lemma:l-kerpr2}, we know that (see \eqref{eq:explicit-lxyab-2}),
$$
		l^{y,x}_{\bf{a},\bf{b}}  =(-1)^{k+l}[\omega^{b_l-1}_{12}\mid\omega_{23}\mid\dots \mid\omega^{b_1-1}_{12}\mid\omega_{23}\mid\omega^{a_k-1}_{12}\mid\omega_{24}\mid\dots \mid\omega^{a_1-1}_{12}\mid\omega_{24}]
$$
and that the depth of $\bf{b}$ corresponds to the number of $\omega_{23}$'s when one restricts $l^{y,x}_{\bf{a},\bf{b}}$ to the space $\ker \pr_2$. Dually, $x_{23}$ is a grading of $\ker \pr_2$.
\end{proof}

  Note that for $\gamma \in \ker \pr_2$, the images $\pi_0^{2,3}(\gamma)$ and $(\pi_1^p)^{2,3}(\gamma)$ are identified, under the map \eqref{eq:exact sequence}, with $\gamma^{(0)}$ and $\gamma^{(1)}$ respectively, where $\gamma = \sum_{i \in \mathbb{N}} \gamma^{(i)}$ is its decomposition into $x_{23}$-graded components. 

\begin{lemma} \label{lemma: the quotient}
	For any $\gamma\in \ker \pr_{2}\subset U\mathfrak{p}^{2,3}_5 $, 
	$
		l^{y,x}_{{\bf{a}},(b_1)}\left( (\pi^{2,3}_0+(\pi^p_1)^{2,3})(\gamma) \right)=l^{y,x}_{{\bf{a}},(b_1)}(\gamma).
	$
\end{lemma}

\begin{proof}
 $\gamma$ decomposes as the sum of its graded components,
$
    \gamma=\pi^{2,3}_0(\gamma)+(\pi_1^p)^{2,3}(\gamma) +\sum_{i \geq 2} \gamma^{(i)}, \text{ where } \gamma^{(i)}\in (I^{2,3})^2.
$
The result follows from Lemma \ref{lemma:filtration}.
\end{proof}

\begin{lemma} \label{lemma:exact_pair}
If $l^{y,x}_{\bold a,(b_1)}(\alpha^{(1)})=0$ for all indices $\bold a, (b_1)$, then $\alpha^{(1)}=0$.
\end{lemma}

\begin{proof}
Recall that $\alpha^{(1)}=(\pi^p_1)^{2,3}(\alpha)$. By Lemma \ref{lemma:gamma-1-P}, $\alpha^{(1)}=P(\ad_{x_{12}}, \ad_{x_{24}})x_{23}$ for some $P \in k\langle \langle x_0, x_1 \rangle \rangle$ and by \eqref{eq:alpha-graded-component-1}, $\alpha^{(1)} \in V^Y_{\kv}$, so every monomial $w$ of $P$ contains at least one occurrence of $x_{24}$: it can thus be written uniquely as $x_{12}^{s_w}x_{24}w'$, where $w'$ is a word in $x_{12}$ and $x_{24}$. 

Assume that $\alpha^{(1)} \neq 0$, then $P \neq 0$. Among the terms with nonvanishing coefficient in $P$, we choose the monomial $w_P = x_{12}^{s_P} x_{24} w_P'$ with $s_P$ minimal among the monomials.

Now, consider the word $x_{12}^{s_P}x_{23}\overline{x_{24}w'_P} = x_{12}^{s_P}x_{23}\overline{w'_P}x_{24}$, where $\overline{w'_P}$ denotes the word $w'_P$ with its letters in reversed order. As this word ends in $x_{24}$, we can associate to it $\textbf{a} = (a_1, \ldots, a_k)$ such that $\overline{w'_P}x_{24} = x_{12}^{a_k-1}x_{24} \ldots x_{12}^{a_1-1}x_{24}$.
Since the action of $\ad$ expands any word of the form $w(\ad_{x_{12}}, \ad_{x_{24}})x_{23}$ into the form $L x_{23} \overline{R}$, where $L,R$ are complementary subwords each read in the order given by $w$, a monomial $w$ of $P$ contributes to this word only if it has $s_P$ letters $x_{12}$ inserted into $x_{24}w'_P$, and it contributes only once. Such a monomial then has prefix $x_{12}^tx_{24}$, so $s_w=t \leq s_P$, which the minimality assumption forces to be an equality: $t=s_P$, hence $w=w_P$. Thus, $w_P$ is the only contributing monomial and the coefficient in $\alpha^{(1)}$ is $(-1)^{\lvert w_P \rvert - s_P} c_{w_P}(P) \neq 0$.
Finally, by \eqref{eq:explicit-lxyab-2},
\begin{align*}
    l^{y,x}_{\bold{a},(s_P+1)}(\alpha^{(1)}) = (-1)^{\wt(\textbf{a})+\dpt(\textbf{a})+1} \underbrace{c_{w_P}(P)}_{\neq 0} \neq 0.
\end{align*}

\end{proof}

\begin{proposition}
\label{prop:one_leg}
If $\psi\in \mathrm{Skew}$ and $c_{x_0x_1}(\psi)=0$, then the following two conditions are equivalent:
\begin{enumerate}[(i)] 
    \item \label{prop:one_leg-1}
    $l^{y,x}_{\bold a,(b_1)}(\alpha)=0$, for all $\bold a, (b_1)$;
    \item \label{prop:one_leg-2}
    $\psi\in \mathfrak{rc}_0$.
\end{enumerate}
\end{proposition}

\begin{proof}
It follows from Lemmas \ref{lemma: the quotient} and \ref{lemma:pi_0} that
$
l^{y,x}_{\bold a,(b_1)}(\alpha)=l^{y,x}_{\bold a,(b_1)}(\pi^{2,3}_0(\alpha)+(\pi_1^p)^{2,3}(\alpha))=l^{y,x}_{\bold a,(b_1)}((\pi^p_1)^{2,3}(\alpha)).
$

The condition \eqref{prop:one_leg-1} implies by Lemma \ref{lemma:exact_pair} that $(\pi^p_1)^{2,3}(\alpha)=0$, which implies $\psi\in \mathfrak{rc}_0$ by Proposition \ref{prop:kernel and rc}. Conversely, \eqref{prop:one_leg-2} implies by Proposition \ref{prop:kernel and rc} that $(\pi^p_1)^{2,3}(\alpha)=0$, which implies $l^{y,x}_{\bold a, (b_1)}(\alpha)=0$ for all indices $\textbf{a},(b_1)$.
\end{proof}

\subsection{Change of variable} \label{subsection:3_2}
Recall that $\mathfrak{p}_4$ is free on the variables $x=x_{12}=x_{34}$ and $y=x_{23}=x_{14}$, and observe that it is also free on the generators $x=x_{12}=x_{34}$ and $z=x_{13}=x_{24}$, as the only relation is $x+y+z=0$. Each choice of generating set thus yields a free Lie algebra; the two we will consider are related by the change of variables illustrated in the diagram below. For $\psi \in \mathfrak{fr}_k(x_0,x_1)$ with $c_{x_0}(\psi)=c_{x_1}(\psi)=0$:
$$
\begin{tikzcd}
	\psi(x_0,x_1)\in\mathfrak{fr}_k(x_0,x_1)\arrow{r}{x_0\mapsto x_{12},x_1\mapsto x_{23}}\arrow{d}{\text{change of variable}}&[3em]\mathfrak{p}_4\arrow{d}{\id}\arrow{r}{\text{coface maps } \eqref{eq:coface-maps-1}}&[3em]\mathfrak{p}_5\arrow{d}{\id}\\
	\psi(-x_0-x_1,x_1)\in\mathfrak{fr}_k( x_0,x_1)\arrow{r}{x_0\mapsto x_{13},x_1\mapsto x_{23}}&\mathfrak{p}_4\arrow{r}{\text{coface maps } \eqref{eq:coface-maps-2}}&\mathfrak{p}_5
\end{tikzcd}.
$$

In the previous Subsection \ref{subsection:3_1}, we studied $\mathfrak{fr}_k(x_{12},x_{23})$, i.e. the top line of this diagram, and in the present subsection, we apply the same method to study the bottom line, $\mathfrak{fr}_k(x_{13},x_{23})$. Most results follow analogously, owing to the dihedral symmetry of the defect $\alpha$.\\

Recall that $U\mathfrak{p}^{3,4}_5:=\ker (\pr_3\circ \pr_4)$ and that $I^{3,4}$ denotes the ideal generated by $x_{34}$ in $U\mathfrak{p}_5^{3,4}$. $U\mathfrak{p}^{3,4}_5$ admits the presentation
$$
	U\mathfrak{p}^{3,4}_5=\langle x_{14},x_{24},x_{34},x_{23},x_{13} \mid R^{3,4} \rangle,
$$
where the $R^{3,4}$ relations are as in Lemma \ref{eq:4T_relations}, with indices $2$ and $4$ exchanged.

\begin{proposition}
	The assignment
	\begin{align*}
		\pi^{3,4} \colon U\mathfrak{p}^{3,4}_5 &\to (A \otimes A \oplus A, \cdot_{\rho_{\rm{KKS}}}) \\
		x_{14} &\mapsto x_0 \otimes 1, \quad x_{24}\mapsto x_1\otimes 1  \\
		x_{13} &\mapsto 1 \otimes x_0,\quad x_{23}\mapsto 1\otimes x_1 \\
		x_{34} &\mapsto -e,
	\end{align*}
	where $e := 0 \oplus 1 \in A \otimes A \oplus M$, extends to an algebra homomorphism. Furthermore, the map $\pi^{3,4}$ factors through the canonical projection $U\mathfrak{p}^{3,4}_5 \to U\mathfrak{p}^{3,4}_5 /(I^{3,4})^2$.
\end{proposition}	

\begin{proof}
    The proof is that of Proposition \ref{prop:pi-2-3}, with the indices $2$ and $3$ replaced by $3$ and $4$.
\end{proof}

By composing $\pi^{3,4}$ with the natural projections 
$A \otimes A \oplus M \to A \otimes A$ and $A \otimes A \oplus M \to  M$, we define the maps
$
\pi^{3,4}_0 \colon U\mathfrak{p}^{3,4}_5 \to A \otimes A, 
$ and 
$
\pi^{3,4}_1 \colon U\mathfrak{p}^{3,4}_5 \to M.
$
Similarly to Proposition \ref{prop:projection}, the map $\pi^{3,4}_0$ descends to an isomorphism 
$$
        \pi^{3,4}_0:U\mathfrak{p}_5^{3,4}/I^{3,4} \overset{\simeq}{\longrightarrow} A\otimes A
$$
and on $I^{3,4}$ the map $\pi^{3,4}_1=F^{3,4}\circ (\pi^p_1)^{3,4}$ where $(\pi^p_1)^{3,4}:I^{3,4}\to I^{3,4}/(I^{3,4})^2$ is the projection map and $F^{3,4}~:~I^{3,4}~/~(I^{3,4})^2~\to~M$ is the bimodule map determined by $F^{3,4}(x_{34})=-1.$ Using a split short exact sequence of bimodules analogous to \eqref{eq:exact sequence}, we can write
\begin{equation*}
\begin{split}
    	U\mathfrak{p}^{3,4}_5/(I^{3,4})^2 \simeq A\otimes A\oplus I^{3,4}/(I^{3,4})^2\\
        a\mapsto (\pi^{3,4}_0(a),(\pi^p_1)^{3,4}(a)).
\end{split}
\end{equation*}

We now introduce another defect $\widehat{\alpha}$ of the pentagon equation of  $\eta(x_0,x_1):=\psi(-x_0-x_1,x_1) \in~\mathfrak{fr}_k(x_0,x_1)$:
$$
 \alpha(\psi) = \widehat{\alpha}(\eta(x_0,x_1))=\eta(x_{13},x_{23})+\eta(x_{14},x_{24}+x_{34})+\eta(x_{24},x_{34})
-\eta(x_{14}+x_{24},x_{34})-\eta(x_{13}+x_{14},x_{23}+x_{24}),
$$
which is derived from $\alpha(\psi)$ in a similar way to \eqref{eq:defect-skew-2}, using the skew symmetry of $\psi$ and additional identities $x_{45}=-x_{14}-x_{24}-x_{34}$ and $-x_{13}-x_{14}=x_{12}+x_{15}$ of $\mathfrak{p}_5$.
$\widehat{\alpha}(\eta)$ is a linear combination of the following coface maps from $k\langle \langle x_{13},x_{23} \rangle \rangle$ to $U\mathfrak{p}_5$:
\begin{equation}
\label{eq:coface-maps-2}
\begin{split}
	c_{1,2,34}:\quad &\eta(x_{13},x_{23})\mapsto \eta(x_{13}+x_{14},x_{23}+x_{24})\\
	c_{2,3,4}:\quad &\eta(x_{13},x_{23})\mapsto \eta(x_{24},x_{34})\\
	c_{12,3,4}:\quad&\eta(x_{13},x_{23})\mapsto \eta(x_{14}+x_{24},x_{34})\\
	c_{1,2,3}:\quad &\eta(x_{13},x_{23})\mapsto \eta(x_{13},x_{23})\\
	c_{1,23,4}:\quad&\eta(x_{13},x_{23})\mapsto \eta(x_{14}, x_{24}+x_{34}).
\end{split}
\end{equation}

\begin{lemma}\label{lemma:cabling34}
	Let $\eta$ be a Lie series, then the following relations hold: 
	\begin{align*}
		&\pi^{3,4}_1\circ c_{1,2,34}(\eta(x_0,x_1))=\mu(\eta(x_0,x_1));\\
        & \pi^{3,4}_1\circ c_{2,3,4}(\eta(x_0,x_1))=-(d^R_1\eta)(x_1,0);\\
		&\pi^{3,4}_1\circ c_{12,3,4}(\eta(x_0,x_1))=-(d^R_1\eta)(x_0+x_1,0);\\
		&\pi^{3,4}_1\circ c_{1,2,3}(\eta(x_0,x_1))=0;\\
        & \pi^{3,4}_1\circ c_{1,23,4}(\eta(x_0,x_1))=-d^R_1(\eta).
	\end{align*}
\end{lemma}

\begin{proof}
	The first and the last equalities are proved in the same way as Propositions $5.8$ and $5.10$ in \cite{ANR}. We present the proof of the second equality; the remaining ones are analogous and rely on the $A\otimes A$ bimodule structure described in \eqref{eq:bimodule}. It holds that $\pi^{3,4}_1\circ c_{2,3,4}(\eta(x_0,x_1))=\pi^{3,4}_1(\eta(x_{24},x_{34}))$ and the map $\pi^{3,4}_1=F^{3,4}\circ(\pi^p_1)^{3,4}$ factors through the quotient  $I^{3,4}/(I^{3,4})^2$, so it suffices to look at the terms which only contain one $x_{34}$. From the bimodule structure, $\pi_0^{3,4}(x_{24})=x_1 \otimes 1$  acts from left by $\varepsilon(x_{24})=0$ and from right by multiplication by $x_{24}$. Hence, only the terms with a single $x_{34}$ term and beginning with $x_{34}$ can contribute. These terms are exactly $-(d^R_1\eta)(x_1,0)$. 
\end{proof}

As a direct corollary, we get the following results.

\begin{proposition}
If $\pi^{3,4}_1\left( \widehat{\alpha}(\eta(x_0,x_1))\right)=0$, then 
$$
\mu(\eta(x_0,x_1))=d^R_1(\eta(x_0,x_1))(x_0+x_1,0)-d^R_1(\eta(x_0,x_1)) -d^R_1(\eta(x_0,x_1))(x_1,0).
$$
\end{proposition}


\begin{corollary}\label{cor:equation34}
If $\pi^{3,4}_1\left(\widehat{\alpha}(\psi(-x_0-x_1,x_1))\right)=0$, then
$$
		\mu(\psi(-x_0-x_1,x_1)) =(d^R_1 \psi(-x_0-x_1,x_1))(x_0+x_1,0) -d^R_{1}(\psi(-x_0-x_1,x_1))
        -d^R_1(\psi(-x_0-x_1,x_1))(x_1,0).
$$
\end{corollary}

Let $\widehat{\gamma}$ be a Lie series in the ideal $(x_{34})\subset \ker \pr_4\subset U\mathfrak{p}^{3,4}_5$. As $x_{34}$ is a grading of the vector space $\ker \pr_4=\mathfrak{fr}_k (x_{14},x_{24},x_{34})$, we write $\widehat{\gamma}=\sum_{i\in \N}\widehat{\gamma}^{(i)}$, where $\widehat{\gamma}^{(i)}$ is the $i$-th graded component with respect to the $x_{34}$ grading.

\begin{lemma} \label{lemma:polynomial_expression}
Let $\widehat{\gamma}$ be as above, then
$   \widehat{\gamma}^{(1)}=P(\ad_{x_{14}},\ad_{x_{24}})x_{34}
$
for some $P\in k\langle\langle x_0,x_1\rangle\rangle$.
\end{lemma}

\begin{lemma}[Change of variable] \label{lemma:change-var-P}
Using the generating set $\{x_{24}, x_{34}, x_{45}=-x_{14}-x_{24}-x_{34}\}$ of $\ker \rm{pr}_4$,
$
\widehat{\gamma}^{(1)}~=~P'(\ad_{x_{45}},\ad_{x_{24}})x_{34}
$
for some $P'\in k\langle\langle x_0,x_1\rangle\rangle$. 
\end{lemma}

We denote by $\widehat{V}_{\rm{kv}}$ the vector space spanned by the noncommutative polynomials $P(\text{ad}_{x_{45}},\text{ad}_{x_{24}})x_{34}$, where $ P(x_0,x_1)\in k\langle\langle x_0,x_1\rangle\rangle$. We denote by $\widehat{V}_{\mathrm{kv}}^Y$ the vector subspace spanned by those $P(x_0,x_1)$ in which every monomial contains at least one occurrence of the second variable $x_1$.

\begin{lemma}\label{lemma:F_injective34}
	$\eval{F^{3,4}}_{\widehat{V}_{\rm{kv}}}$ is injective.
\end{lemma}
\begin{proof}
	$F^{3,4}$ maps $P(\ad_{x_{45}},\ad_{x_{24}})x_{34}$ to  $-\overline{P}(x_0+x_1,-x_1)$, where $\overline{P}$ denotes $P$ with each of its monomials read from right to left.
\end{proof}

\begin{lemma} 
Let $\textbf{a}=(a_1, \ldots, a_k),\textbf{b}=(b_1, \ldots, b_l)$, then
$
		\eval{l^{x,y}_{\bf{a},\bf{b}}}_{\ker \pr_4} \left( (I^{3,4})^t\cap \ker \pr_4 \right)=0,\quad \text{for } t> \dpt(\bf{b}).
$
\end{lemma}
\begin{proof}
By Lemma \ref{lemma:l-kerpr2}, we know that (see \eqref{eq:explicit-lxyab}), 
\[
\eval{l^{x,y}_{\bf{a},\bf{b}}}_{\ker \pr_4}=(-1)^{k+l}[\omega^{b_l-1}_{45} \mid \omega_{34} \mid \dots \mid \omega^{b_1-1}_{45} \mid \omega_{34} \mid \omega^{a_k-1}_{45} \mid \omega_{24} \mid \dots \mid \omega^{a_1-1}_{45} \mid \omega_{24}].
\]
The number of $\omega_{34}$'s is exactly $\dpt(\textbf{b})$. Dually, $x_{34}$ is a grading of the vector space $\ker \pr_4$.
\end{proof}

Again, for $\widehat{\gamma} \in \ker \pr_4$, the images $\pi_0^{3,4}(\widehat{\gamma})$ and $(\pi_1^p)^{3,4}(\widehat{\gamma})$ are identified with $\widehat{\gamma}^{(0)}$ and $\widehat{\gamma}^{(1)}$ respectively, where $\widehat{\gamma} = \sum_{i \in \mathbb{N}} \widehat{\gamma}^{(i)}$ is its decomposition into $x_{34}$-graded components.

\begin{lemma}\label{lemma: the quotient_hat}
	For any $\widehat{\gamma}\in \ker \pr_4 \subset U\mathfrak{p}^{3,4}_5$,
$
		l^{x,y}_{{\bf{a}},(b_1)}(\widehat{\gamma})=l^{x,y}_{{\bf{a}},(b_1)}\left( (\pi^{3,4}_0+(\pi^p_1)^{3,4})(\widehat{\gamma}) \right).
$
\end{lemma}

\begin{proof}
    The proof is that of Lemma \ref{lemma: the quotient}, with the indices $2$ and $3$ replaced by $3$ and $4$.  
\end{proof}

\begin{lemma}
	\label{lemma: the kernel_hat}
	If $l^{x,y}_{{\bf{a}},(b_1)}(\widehat{\alpha}^{(1)})=0$ for all indices ${\bf{a}},(b_1)$, then $\widehat{\alpha}^{(1)}=0$.
\end{lemma}

\begin{proof}
The proof is analogous to that of Lemma \ref{lemma:exact_pair}, replacing the variables $x_{12}$ and $x_{23}$ by $x_{45}$ and $x_{34}$ respectively ($x_{24}$ is unchanged), using Lemma \ref{lemma:change-var-P} to write $\widehat{\alpha}^{(1)} = P'(\ad_{x_{45}},\ad_{x_{24}})x_{34}$,
for some $P' \in k\langle \langle x_0, x_1 \rangle \rangle$. By Lemma \ref{lemma:skew}\eqref{lemma:skew-ii}, $\widehat{\alpha}(\psi(x_\infty,x_1)) = \alpha(\psi) \in (x_{24})$, so every monomial of $P'$ contains at least one occurrence of $x_{24}$ and $\widehat{\alpha}^{(1)} \in \widehat{V}^Y_{\kv}$. Finally, we conclude using equation \eqref{eq:explicit-lxyab}.
\end{proof}

\begin{proposition}
If $\psi\in \mathrm{Skew}$, then the following two conditions are equivalent:
	\begin{enumerate}[(i)] 
		\item \label{prop:lab-coaction-1}
        $l_{{\bf{a}},(b_1)}^{x,y}(\alpha)=0,\quad \text{for all } \textbf{a},(b_1)$;
		\item \label{prop:lab-coaction-2}
        $\begin{aligned}[t]
		\mu(\psi(-x_0-x_1,x_1)) & =d^R_1(\psi(-x_0-x_1,x_1))(x_0+x_1,0)\\ & -d^R_{1}(\psi(-x_0-x_1,x_1))-d^R_1(\psi(-x_0-x_1,x_1))(x_1,0).
	\end{aligned}$
	\end{enumerate}
\end{proposition}

\begin{proof}
We know that 
$
l^{x,y}_{{\bf{a}},(b_1)}(\alpha(\psi(x_0,x_1)))=l^{x,y}_{{\bf{a}}, (b_1)}(\widehat{\alpha}(\psi(-x_0-x_1,x_1))),
$
and by skew symmetry, $\pi^{3,4}_0(\widehat{\alpha})=0$. It follows from Lemma \ref{lemma: the quotient_hat} that
\begin{align*}
l^{x,y}_{{\bf{a}}, (b_1)}(\widehat{\alpha}(\psi(-x_0-x_1,x_1))
=l^{x,y}_{{\bf{a}},(b_1)}(\pi^{3,4}_0(\widehat{\alpha})+(\pi^p_1)^{3,4}(\widehat{\alpha}))
=l^{x,y}_{{\bf{a}},(b_1)}((\pi^p_1)^{3,4}(\widehat{\alpha})).
\end{align*}

From \eqref{prop:lab-coaction-1} to \eqref{prop:lab-coaction-2}, we have by Lemma \ref{lemma: the kernel_hat} that $(\pi^p_1)^{3,4}(\widehat{\alpha})=0$, which implies $\pi^{3,4}_1(\widehat{\alpha})=F^{3,4}\circ (\pi^p_1)^{3,4}(\widehat{\alpha})=0$. The result \eqref{prop:lab-coaction-2} then follows from Corollary \ref{cor:equation34}.

From \eqref{prop:lab-coaction-2} to \eqref{prop:lab-coaction-1}, the condition \eqref{prop:lab-coaction-2} is equivalent to $\pi^{3,4}_1(\widehat{\alpha})=0$, as $F^{3,4}$ is injective in the subspace $\widehat{V}_{\kv}$ by Lemma \ref{lemma:F_injective34}. This implies $(\pi^p_1)^{3,4}(\widehat{\alpha})=0$, and then $l^{x,y}_{{\bf{a}},(b_1)}(\alpha)=0$.
\end{proof}

\subsection{Dihedral symmetry and equivalence of the formulas} \label{subsection:3_3}

In this subsection, we consider the defect of the pentagon equation of a skew-symmetric Lie series $\psi$ in the form 
\begin{equation*}
\alpha=\psi(x_{12},x_{23})+\psi(x_{23},x_{34})+\psi(x_{34},x_{45})+\psi(x_{45},x_{51})+\psi(x_{51},x_{12}).
\end{equation*}

Recall from Subsection \ref{subsubsec:skew} 
that it follows from the action of the dihedral group on $U\mathfrak{p}_5$ and the skew symmetry of $\psi$ that
$ 
\alpha^{\sigma} = \alpha,$
and $\alpha^{\tau} = -\alpha$ (see \eqref{eq:dihedral_symmetry}),
where $\tau$ and $\sigma$ are the generators of the dihedral group $D_5$ defined in \eqref{eq:dihedral_generator}.

\begin{lemma}\label{lem:dihedral_forms}
	For any element $\varphi\in U\mathfrak{p}_5$ and index $\textbf{a},\textbf{b}$,
$
		l^{y,x}_{\bf{a},\bf{b}}(\varphi)=l^{x,y}_{\bf{a},\bf{b}}(\varphi^{\tau}).
$
\end{lemma}
\begin{proof}
    The formula follows from the graded duality between the bar construction and $U\mathfrak{p}_{r+3}$. By duality, the morphism $\tau$ induces a map on $\mathcal{V}(\mathcal{M}_{0,5})$, which maps $l^{x,y}_{\bold a,\bold b}$ to $l^{y,x}_{\bold a,\bold b}$.  
\end{proof}

\begin{proposition}
If $\psi\in \mathrm{Skew}$, 
then the following are equivalent:
\begin{enumerate}[(i)] 
    \item $l^{y,x}_{\bold a,(b_1)}(\alpha)=0$, \quad for all $\bold a, (b_1)$;
    \item $l^{x,y}_{\bold a,(b_1)}(\alpha)=0$, \quad for all $\bold a, (b_1).$
\end{enumerate}
\end{proposition}

\begin{proof}
It follows from Lemma \ref{lem:dihedral_forms} and the dihedral symmetry of $\alpha$ 
that
$
l^{y,x}_{\bold a,(b_1)}(\alpha)=l^{x,y}_{\bold a, (b_1)}(\alpha^{\tau})=-l^{x,y}_{\bold a,(b_1)}(\alpha).
$
\end{proof}

\section{Relation between $\mathfrak{dmr}_0$ and $\mathfrak{rc}_0$}
\label{sec:linking-dmr-rc}

In this section, we study in greater detail the polylogarithms $l^{y,x}_{\bf{a}, \bf{b}}$ which have a second index \textbf{b} of depth $1$, and leverage their properties 
to investigate the relation between the Lie algebras $\mathfrak{rc}_0$ and $\mathfrak{dmr}_0$. For convenience, we use the shorthand $l^{y,x}_{k,l}:=l^{y,x}_{(1,\dots,1),(1,\dots,1)}$ to denote the word that has $k$ $1$'s in $\bf{a}$ and $l$ $1$'s in $\bf{b}$. 

\subsection{Polylogarithm calculations with $\dpt({\textbf{b}})=1$}

\begin{lemma}\label{lemma:451123_1}
Let $\psi\in \mathfrak{dmr}_0$, then
$
l^{y,x}_{k,1}(\psi_{451}+\psi_{123})=(-1)^{k+1}c_{x^k_0x_1}(\psi).
$
\end{lemma}

\begin{proof}
By Proposition \ref{prop:nonadmissible1111}, applied with \textbf{a}$=(1)$ and \textbf{b}$=\underbrace{(1, \ldots, 1)}_k$, we have
\begin{align*}
   l^{y,x}_{\underbrace{(1,\dots,1)}_{k},(1)}(\psi_{451}+\psi_{123})-l_{ \underbrace{(1,\dots,1)}_{k+1}}(\psi)=(-1)^{k+1}c_{x^k_0x_1}(\psi), 
\end{align*}
since $l_{(1, \ldots, 1)}(\psi)=0$ for a Lie series.
\end{proof}

\begin{lemma}\label{lemma:432_1}
The parts of $l^{y,x}_{k,l}$ that consist only of $\omega_{23}$ and $\omega_{34}$ satisfy the following recursion relations
	\begin{enumerate}
		\item \label{item:formula1} If $k>1,l>1$, $l^{y,x}_{k,l}=-[\omega_{34}\mid l^{y,x}_{k-1,l}]+[\omega_{34}\mid l^{y,x}_{k,l-1}]-[\omega_{23}\mid l^{y,x}_{k,l-1}]$;
        
		\item \label{item:formula2}
        If $k\ge 1,l=1$, $l^{y,x}_{k,1}=(-1)^{k+1}[\omega^k_{34}\mid \omega_{23}]$;
        
		\item If $k=1,l>1$, $l^{y,x}_{1,l}=(-1)^{l+1}[\omega_{34}\mid \omega^l_{23}]+[\omega_{34}\mid l^{y,x}_{1,l-1}]-[\omega_{23}\mid l^{y,x}_{1,l-1}]$.
	\end{enumerate}
    Moreover, $l^{y,x}_{k,1}(\psi_{432})=(-1)^{k+1}c_{x^k_0 x_1}(\psi)$.
\end{lemma}

\begin{proof}
To determine an expression for $l_{k,l}^{y,x}$, we exploit the differential equations given in \eqref{equation:2_variables_different_equation}, switching the roles of $x$ and $y$. Since we are interested only in the parts consisting of $\omega_{23}$ and $\omega_{34}$, it suffices to keep the terms in $\frac{dy}{1-y}$ and $\frac{dx}{1-x}$ in these differential equations. We detail this procedure below for the case \ref{item:formula1}, the other two cases being similar. Let $Li_{k,l}(x,y)$ be shorthand for $Li_{\underbrace{(1, \ldots,1)}_k,\underbrace{(1, \ldots,1)}_l}(x,y)$.

In the first case, $a_k =1, k>1, b_l=1, l>1$, so \eqref{equation:2_variables_different_equation} tells us 
\begin{equation*}
\begin{split}
    \frac{d}{dy} Li_{k,l}(y,x) & = \frac{1}{1-y}Li_{k-1,l}-(\frac{1}{y}+\frac{1}{1-y})Li_{k,l-1}(y,x),\\
    \frac{d}{dx} Li_{k,l}(y,x) & = \frac{1}{1-x}Li_{k,l-1}(y,x),
\end{split}   
\end{equation*}
which translates to, keeping only the terms in $\frac{dy}{1-y}$ and $\frac{dx}{1-x}$ and using the notation above,
\begin{equation*}
    l_{k,l}^{y,x} = [-\omega_{34} \mid l_{k-1,l}^{y,x}] + [\omega_{34} \mid l_{k,l-1}^{y,x}] - [\omega_{23} \mid l_{k,l-1}^{y,x}],
\end{equation*}
as required.\\

Finally, we compute $l^{y,x}_{k,1}(\psi_{432})=\langle (-1)^{k+1}[\omega^k_{34}\mid \omega_{23}], \psi_{432}\rangle=(-1)^{k+1}c_{x^k_0x_1}(\psi).$
\end{proof}

\begin{proposition}\label{prop:defect nonadmissible}
	If $\psi\in \mathfrak{dmr}_0$, then
	$
		l^{y,x}_{{\bf{a}},(b_1)}(\alpha)=0, \text{ for any ${\bf{a}},(b_1)$}.
	$
\end{proposition}

\begin{proof}
The case $\bold a, (b_1)$ $\ne (1,\ldots,1),(1)$ follows from Proposition \ref{prop:defect_admissible}. For the case $\bold a, (b_1)$ $=(1,\ldots,1),(1)$, we have by Lemmas \ref{lemma:432_1} and \ref{lemma:451123_1},
\begin{equation*}
l^{y,x}_{k,1}(\psi_{451}+\psi_{123})=l^{y,x}_{k,1}(\psi_{432})=(-1)^{k+1}c_{x^k_0x_1}(\psi) \quad \text{ and } \quad l^{y,x}_{k,1}(\psi_{215})=l^{y,x}_{k,1}(\psi_{543})=0,
\end{equation*}
therefore
$
l^{y,x}_{k,1}(\alpha)=l^{y,x}_{k,1}(\psi_{451}+\psi_{123}-\psi_{432}-\psi_{215}-\psi_{543})=0.
$
\end{proof}

\begin{proposition} \label{prop:one_loop_double_shuffle}
    If $\psi\in \mathfrak{fr}_k(x_0,x_1)$ with $c_{x_0}(\psi)=c_{x_1}(\psi)=0$, then the following two conditions are equivalent:
	\begin{enumerate}[(i)] 
		\item \label{prop:one_loop_double_shuffle-1}
        $l^{y,x}_{{\bf{a}},(b_1)}(\alpha)=0,\quad \text{for any ${\bf{a}},(b_1)$}$;
		\item \label{prop:one_loop_double_shuffle-2}
        For any ${\bf{a}},(b_1)$, 
        $
            \sum_{\sigma\in \text{Sh}^{\le (1,k)}}
        l_{\sigma\left((b_1),{\bf{a}}\right)}(\psi_{*})=0 \text{ and } \sum_{\sigma\in \text{Sh}^{\le (k,1)}}l_{\sigma\left({\bf{a}},(b_1)\right)}(\psi_{*})=0. 
        $
	\end{enumerate}
\end{proposition}

\begin{proof}
To show \eqref{prop:one_loop_double_shuffle-1} implies \eqref{prop:one_loop_double_shuffle-2}, we distinguish two cases. In the case ${\bf{a}},(b_1)\ne (1,\dots,1),(1)$, the condition \eqref{prop:one_loop_double_shuffle-1} is equivalent to 
$$
		l^{y,x}_{{\bf{a}},(b_1)}(\psi_{451}+\psi_{123})=l^{y,x}_{{\bf{a}},(b_1)}(\psi_{215})\overset{ \text{Lemma } \ref{lemma:polylogs-compilation}}{=}l_{{\bf{a}}b_1}(\psi),\quad \text{for }{{\bf{a}}b_1}\ne (1,\dots,1).
$$
We then evaluate equation \eqref{eq:stuffle 2 variable} in $\psi_{451}+\psi_{123}$, which yields
$$
\sum_{\sigma\in \text{Sh}^{\le (1,k)}}l_{\sigma\left((b_1),{\bf{a}}\right)}(\psi)=\sum_{\sigma\in \text{Sh}^{\le (1,k)}}l_{\sigma\left((b_1),{\bf{a}}\right)}(\psi_{*})=0, \quad \text{for }{{\bf{a}}b_1}\ne (1,\dots,1).
$$

For the other case, when $\textbf{a},(b_1)=(1, \ldots, 1), (1)$, the condition \eqref{prop:one_loop_double_shuffle-1} is equivalent to
$$
l^{y,x}_{k,1}(\psi_{451}+\psi_{123})-l^{y,x}_{k,1}(\psi_{215})=l^{y,x}_{k,1}(\psi_{432})
$$
which, by Lemma \ref{lemma:432_1} is the same as
$$
l^{y,x}_{(1,\dots,1),(1)}(\psi_{451}+\psi_{123}) =l_{\underbrace{(1,\dots,1)}_{k+1}}(\psi)+(-1)^{k+1}c_{x_0^{k}x_1}(\psi).
$$

We then evaluate equation \eqref{eq:stuffle 2 variable} in $\psi_{451}+\psi_{123}$, to get
$$
\sum_{\sigma\in \text{Sh}^{\le (1,k)}}l_{\sigma\left((b_1),{\bf{a}}\right)}(\psi)+(-1)^{k+1}c_{x^k_0x_1}(\psi)=\sum_{\sigma\in \text{Sh}^{\le (1,k)}}l_{\sigma\left((b_1),{\bf{a}}\right)}(\psi_{*}).
$$
By \eqref{eq:stuffle 2 variable}, this implies the stuffle relations of $\psi_{*}$ for such indices ${\bf{a}},(b_1)$. Since the stuffle product is commutative, the relations
$$
    \sum_{\sigma\in \text{Sh}^{\le (1,k)}}l_{\sigma\left((b_1),{\bf{a}}\right)}(\psi_{*})=0 \quad \text{and} \quad \sum_{\sigma\in \text{Sh}^{\le (k,1)}}l_{\sigma\left({\bf{a}},(b_1)\right)}(\psi_{*})=0
$$
are the same. \\

The converse is proved by inverting the above step.
\end{proof}

Recall that
the \textit{meta-abelian quotient} of $\psi$ is defined by
$
	B_{\psi}(x_0,x_1):=\left(\psi{x_1}x_1\right)^{\text{ab}},
$
where $h \mapsto h^{\text{ab}}$ is the abelianization map from $k\langle\langle x_0,x_1\rangle\rangle$ to $k[[x_0,x_1]]$. We also introduce the vector space
$$
	\mathfrak{B}:=\{\beta\in k[[x_0,x_1]] \mid \beta(x_0,x_1)=\gamma(x_0)+\gamma(x_1)-\gamma(x_0+x_1),\gamma(s)\in s^2k[[s]]\}.
$$

\begin{proposition}\label{prop:main_property}
Let $\psi\in \mathfrak{fr}_k(x_0,x_1)$ and $c_{x_0}(\psi)=c_{x_1}(\psi)=0$. If $l^{y,x}_{{\bf{a}},(b_1)}(\alpha)=0$ for any ${\bf{a}},(b_1)$, then 

\begin{enumerate}[(i)]
    \item \label{eq:odd_coefficent}
    $c_{x_0^{n+1}x_1}(\psi)=0, \quad \text{for } n \text{ even}, n\ge 2; $
    \item \label{eq:abelian quotient}
	$B_{\psi}(x_0,x_1)\in \mathfrak{B}$.
\end{enumerate}

\end{proposition}

\begin{rmk}
 The proof of \eqref{eq:odd_coefficent} is the same as Racinet's \cite{Racinet2002}, and the proof of \eqref{eq:abelian quotient} is the same as Furusho's \cite{Furusho2011} for $\mathfrak{dmr}_0$.
The main purpose of repeating them here is to emphasise that the condition $l^{y,x}_{{\bf{a}},(b_1)}(\alpha)=0$ is enough for those two properties, thanks to Proposition \ref{prop:one_loop_double_shuffle}.
\end{rmk}

\begin{proof}
We start by proving point \eqref{eq:odd_coefficent}.	Since $\psi$ is a Lie series, the antipode satisfies 
\begin{align*}
    S(\psi)=-\psi \quad \text{ and } \quad S(x_1x^{n-2}_0x_1)=(-1)^{{n}}x_1x^{n-2}_0x_1.
\end{align*}
We have the equation
\begin{equation}\label{equation: proof_odd}
	c_{x_1x^{n-2}_0x_1}(\psi)+(-1)^{n}c_{x_1x^{n-2}_0x_1}(\psi)=0,
\end{equation}
as well as the shuffle product relation
\begin{align*}
	x_1\shuffle x^{n-2}_0x_1=x^{n-2}_0x_1x_1+\sum_{p=1}^{n-1}x^{p-1}_0x_1x^{n-1-p}_0x_1,
\end{align*}
and since $\psi$ is a Lie series, we know that
\begin{equation}\label{eq:even_middle_1}
	c_{x^{n-2}_0x_1x_1}(\psi)+\sum_{p=1}^{n-1}c_{x^{p-1}_0x_1x^{n-1-p}_0x_1}(\psi)=0.
\end{equation}

We have the following stuffle relation,
\begin{equation}\label{eq:even_middle_2}
	x^{p-1}_0x_1\star x^{n-p-1}_0x_1=x^{n-1}_0x_1+x^{p-1}_0x_1 x^{n-p-1}_0x_1+x^{n-p-1}_0x_1x^{p-1}_0x_1
\end{equation}
and by Proposition \ref{prop:one_loop_double_shuffle}, $\psi$ dually satisfies the relation
\begin{equation}\label{eq:even_middle_2_2}
	c_{x^{n-1}_0x_1}(\psi)+c_{x^{p-1}_0x_1 x^{n-p-1}_0x_1}(\psi)+c_{x^{n-p-1}_0x_1x^{p-1}_0x_1}(\psi)=0.
\end{equation}
Combining the relations \eqref{eq:even_middle_1} and  \eqref{eq:even_middle_2_2} for $p=1,\dots,n-1$, we get
$
	2c_{x^{n-2}_0x_1x_1}(\psi)=(n-1)c_{x^{n-1}_0x_1}(\psi)
$
and using \eqref{eq:even_middle_2} for $p=1$, we get $2c_{x_1x^{n-2}_0x_1}(\psi)=-(n+1)c_{x^{n-1}_0x_1}(\psi)$.
Finally, by equation \eqref{equation: proof_odd},
$$
	c_{x^{n-1}_0x_1}(\psi)+(-1)^nc_{x^{n-1}_0x_1}(\psi)=0,\quad \text{if}\quad n\ge 3.
$$

We now prove point \eqref{eq:abelian quotient}. By Proposition \ref{prop:one_loop_double_shuffle}, we know that $\psi_{*}$ satisfies the stuffle relation of type $({\bf{a}},b_1)$ and $(b_1,{\bf{a}})$. Summing up all pairs $(\bf{a},\bf{b})$ satisfying $\dep({\bf{a}})=1$ and $\text{wt}({\bf{a}})+\text{wt}({\bf{b}})=w$, we get
$$
	\sum_{\substack{\text{wt}({\bf{a}})=w,\\\text{dp}({\bf{a}})=k+1}}(k+1)l_{\bf{a}}(\psi_{*})+\sum_{\substack{\text{wt}({\bf{a}})=w,\\ \text{dp}({\bf{a}})=k}}(w-k)l_{\bf{a}}(\psi_{*})=0,
$$
where $k$ is the weight of $\bold a$. By induction on $k$, we get the relation
$$
	\sum_{\substack{\text{wt}(\bf{a})=w,\\ \text{dp}(\bf{a})=m}}l_{\bf{a}}(\psi)=\begin{cases}
		(-1)^{m-1}\begin{pmatrix}
			w\\m
		\end{pmatrix}\frac{l_{w}(\psi)}{w}, & \quad \text{for \quad $m<w$}\\
		0, &\quad \text{for $m=w$},
	\end{cases}
$$
which proves the property \eqref{eq:abelian quotient}; see \cite[Lemma 5.3]{Furusho2011}.
\end{proof}

\begin{corollary} \label{prop:c-b-rc0}
    If $\psi\in \mathfrak{rc}_0$, then
    \begin{enumerate}[(i)]
        \item \label{eq:even} $c_{x_0^{n+1}x_1}(\psi)=0, \quad \text{for}\quad n\ge 0\quad\text{even};$
        \item \label{eq:abelian_quotient_i} $B_{\psi}(x_0,x_1)\in \mathfrak{B}.$
    \end{enumerate}
\end{corollary}

\begin{proof}
Let $\psi\in \mathfrak{rc}_0$, then Theorem \ref{th:poly_rc_0} implies that $l^{y,x}_{{\bf{a}},(b_1)}(\alpha)=0$ for any ${\bf{a}}, (b_1)$. Then, \eqref{eq:even} for $n\ge 1$ and \eqref{eq:abelian_quotient_i} follow from Proposition \ref{prop:main_property}. The case $n=0$ of \eqref{eq:even} follows from $c_{x_0x_1}(\psi)=0$, by definition of $\mathfrak{rc}_0$. 
\end{proof}

\begin{theorem} \label{theorem:property}
$\mathfrak{rc}_0$ is a Lie algebra with respect to the Ihara bracket \eqref{eq:Ihara_bracket}. Moreover, for any $\psi_1,\psi_2\in \mathfrak{rc}_0$,
\begin{align} \label{eq:Ihara}
	\mu(\{\psi_1,\psi_2\})=-d_0^L(\{\psi_1,\psi_2\})-d_1^R(\{\psi_1,\psi_2\}).
\end{align}
\end{theorem}

\begin{proof}
    Let $\eta \in \rc_0$. By Corollary \ref{prop:c-b-rc0}\eqref{eq:even}, its reduced coaction equation \eqref{reduced_coaction_equation}
\begin{align*}
\mu(\eta)=-r_{\eta}(x_1)+r_{\eta}(-x_0)-d^L_0(\eta)-d^R_1(\eta),
\end{align*}
coincides with
\begin{equation}\label{eq:bar_rc}
		\mu(\eta)=-r_{\eta}(x_1)+r_{\eta}(x_0)-d^L_0(\eta)-d^R_1(\eta).
\end{equation}

It was shown in \cite{Ren2025} that the set of elements $\psi \in \mathrm{Skew}$ satisfying \eqref{eq:bar_rc} forms a Lie algebra  $\overline{\rc_0}$ with respect to the Ihara bracket. It follows that for $\psi_1, \psi_2 \in \rc_0, \psi_1, \psi_2 \in \overline{\rc_0}$, hence $\{\psi_1, \psi_2 \} \in \overline{\rc_0}$. Moreover, $r_{\{\psi_1, \psi_2\}} = 0$, because $\psi_1$ and $\psi_2$ have $x_1$-degree at least $1$, and each component of the Ihara bracket $d_{\psi_1}({\psi_2}), d_{\psi_2}({\psi_1})$ and $[\psi_1, \psi_2]$ has $x_1$-degree at least 2, so $c_{x_0^n x_1}(\{\psi_1, \psi_2\}) = 0$ for all $n$. 
Thus $\{\psi_1, \psi_2\}$ satisfies \eqref{reduced_coaction_equation} and lies in $\rc_0$. \\

Finally, substituting $\{\psi_1, \psi_2\}$ into \eqref{reduced_coaction_equation}, we retrieve \eqref{eq:Ihara} since $r_{\{\psi_1, \psi_2\}} = 0$. 

\end{proof}

\begin{corollary}
    The Lie algebra $\mathfrak{rc}_0$ is isomorphic to $\overline{\mathfrak{rc}_0}$.
\end{corollary}

\begin{proof}
    We deduce from the proof of Theorem \ref{theorem:property} that $\mathfrak{rc}_0$ has the same defining equations as $\overline{\mathfrak{rc}_0}$.
\end{proof}

\subsection{$\mathfrak{rc}_0$ and $\mathfrak{dmr}_0$}
\begin{theorem}[Theorem \ref{theorem:rc_dmr-A}]
\label{theorem:rc_dmr}
	If $\psi\in \mathrm{Skew}$, then the following two conditions are equivalent:
	\begin{enumerate}[(i)]
		\item \label{thm:rco-dmro-equiv-1}
        $\psi\in \mathfrak{dmr}_0$;
		\item \label{thm:rco-dmro-equiv-2}
        $\psi\in \mathfrak{rc}_0$ and for any 
        \textbf{a,b} $\ne (1,\dots,1),(1,\dots,1)$,
		\begin{equation}
        \label{eq:rco-dmro-equiv-2-condition2}
			l^{y,x}_{(a_1,\dots,a_k),(b_1,\dots,b_l)}(\psi_{451}+\psi_{123})=l^{y,x}_{(a_1,\dots,a_k,b_1),(b_2,\dots,b_l)}(\psi_{451}+\psi_{123}).
		\end{equation}
	\end{enumerate}
\end{theorem}

\begin{proof}
We start by proving that \eqref{thm:rco-dmro-equiv-1} implies \eqref{thm:rco-dmro-equiv-2}. Suppose that $\psi \in \mathfrak{dmr}_0 \cap$Skew, then it follows from Proposition \ref{prop:defect nonadmissible} that
$
	   l^{y,x}_{{\bf{a}},(b_1)}(\alpha)=0 \text{ for any } \textbf{a},(b_1), 
$
which implies $\psi\in \mathfrak{rc}_0$ by Proposition \ref{prop:one_leg}. Finally, it follows from Proposition \ref{prop:defect_admissible} and Lemma \ref{lemma:polylogs-compilation} that
$
	l^{y,x}_{\bf{a},\bf{b}}(\psi_{451}+\psi_{123})=l_{\bf{a}\bf{b}}(\psi), \text{ for } ({\bf{a} ,\bf{b}})\ne (1,\dots,1),(1,\dots,1).
$
Hence, we have proven that, for any $(a_1,\dots,a_m,b_1,b_2,\dots,b_n)\ne (1,\dots,1)$,
$$
    l^{y,x}_{(a_1,\dots,a_m),(b_1,\dots,b_n)}(\psi_{451}+\psi_{123}) =l^{y,x}_{(a_1,\dots,a_m,b_1),(b_2,\dots,b_n)}(\psi_{451}+\psi_{123})
	=l_{(a_1,\dots,a_m,b_1,b_2,\dots,b_n)}(\psi).
$$

We now show that \eqref{thm:rco-dmro-equiv-2} implies \eqref{thm:rco-dmro-equiv-1}. Suppose that $\psi\in \mathfrak{rc}_0$ and observe that in view of Theorem \ref{theorem:dmr_defect}, it suffices to consider $\textbf{a},\textbf{b} \neq (1, \ldots, 1),(1, \ldots, 1).$ In that case, it follows from Proposition \ref{prop:one_leg} that
$$
	l^{y,x}_{(a_1,\dots,a_k),(b_1)}(\psi_{451}+\psi_{123})=l_{(a_1,\dots,a_k,b_1)}(\psi)
$$
and together with \eqref{eq:rco-dmro-equiv-2-condition2}, we know that for any $(a_1,\dots,a_k)$,$(b_1,b_2,\dots,b_l)$ $\ne (1,\dots,1)$,$(1,\dots,1)$,
\begin{equation*}
\begin{split}
	l^{y,x}_{(a_1,\dots,a_k),(b_1,\dots,b_l)}(\psi_{451}+\psi_{123}) & = l^{y,x}_{(a_1,\dots,a_k, b_1),(b_2,\dots,b_l)}(\psi_{451}+\psi_{123}) \\
    &  = \ldots = l^{y,x}_{(a_1,\dots,a_k, b_1, \ldots, b_{l-1}),(b_l)}(\psi_{451}+\psi_{123}) 
    = l_{(a_1,\dots,a_k,b_1,b_2,\dots,b_l)}(\psi).
\end{split}
\end{equation*}
 Finally, $\psi\in \mathfrak{dmr}_0$ by Proposition \ref{prop:condition-dmr0}.
\end{proof}

\section{Relations to the Kashiwara--Vergne Problem}
\label{sec:Kashiwara--Vergne}

\subsection{Kashiwara--Vergne Lie algebra}\label{subsection:kv}
Recall that for $A=k\langle\!\langle x_0,x_1\rangle\!\rangle$, we denote by $\operatorname{Der}(A)$ the Lie algebra of derivations of $A$. Define the vector space $\operatorname{tDer}(A):=A\oplus A$, and for $u=(u_0,u_1)\in\operatorname{tDer}(A)$, denote by $\rho(u)\in\operatorname{Der}(A)$
the derivation determined by
\[
\rho(u)(x_0)=[x_0,u_0],
\qquad
\rho(u)(x_1)=[x_1,u_1].
\]
The space $\operatorname{tDer}(A)$ of \textit{tangential derivations} is a Lie algebra, for the Lie bracket given by
$[u,v]=([u,v]_0,[u,v]_1),$ where
\[
[u,v]_i
=
\rho(u)(v_i)-\rho(v)(u_i)+[u_i,v_i],
\qquad i\in\{0,1\}.
\]
With this bracket, the map $\rho:\operatorname{tDer}(A)\to\operatorname{Der}(A)$ is a Lie algebra homomorphism. A tangential derivation
$u=(u_0,u_1)\in\operatorname{tDer}(A)$ is called \textit{special} if $[x_0,u_0]+[x_1,u_1]=0$, or equivalently, $\rho(u)(x_\infty)=0$ for $ x_\infty:=-x_0-x_1$. We denote the Lie subalgebra of \textit{special tangential derivations} by
\[
\operatorname{sDer}(A)
=
\left\{
(u_0,u_1)\in\operatorname{tDer}(A)
\;\middle|\;
[x_0,u_0]+[x_1,u_1]=0
\right\}.
\]
In the same way, we define
$\mathfrak{tder}_2:=\mathfrak{fr}_k(x_0,x_1)\oplus \mathfrak{fr}_k(x_0,x_1)$
and $\mathfrak{sder}_2:=
\left\{
(u_0,u_1)\in\mathfrak{tder}_2
\;\middle|\;
[x_0,u_0]+[x_1,u_1]=0
\right\}.$

Define $|A|:=A/[A,A]$ and consider the projection map $|\cdot|:A\to A/[A,A]$. The \emph{divergence} map is
\begin{align}
\label{eq:divergence}
\rm{div}: \rm{tDer}(A) & \to |A| \nonumber \\
u=(u_0,u_1) & \mapsto |x_0d^R_0(u_0)+x_1d^R_1(u_1)|,
\end{align}
where $d_0^R, d_1^R$ are the right Fox derivatives from Example \ref{example:left_right fox deri}. We now recall the definition of the \emph{Kashiwara--Vergne} Lie algebra $\mathfrak{krv}_2$.

\begin{definition}[\cite{Alekseev2012}] \label{def:krv} $\mathfrak{krv}_2$ consists of the tangential derivations $u\in  \mathfrak{tder}_2 $ which satisfy the following two equations:
\begin{itemize}[leftmargin=5em]
    \item[(krv1)] $u\in \mathfrak{sder}_2$;
    \item[(krv2)] $\text{div}(u)=|f(x_0+x_1)-f(x_0)-f(x_1)|$, for some $f\in k[[x]].$
\end{itemize}
\end{definition}

\subsection{Noncommutative $\rm{krv2}$ equation}

\begin{definition}
    Let $\psi\in \mathfrak{fr}_k(x_0,x_1)$. 
    \begin{enumerate}[(i)]
        \item The \textit{potential function} associated to $\psi$ is defined by
\begin{align*}
h_{\psi}:=x_0\psi(-x_0-x_1,x_0)+x_1\psi(-x_0-x_1,x_1).
\end{align*}

        \item $\psi$ is said to satisfy the \emph{noncommutative krv2 equation} (NCkrv2) if the reduced coaction $\mu$ \eqref{eq:reduced_coaction} of its potential function satisfies the equation
\begin{align}\label{eq:nonKV}
\mu(h_{\psi})=f(x_0+x_1)-f(x_0)-f(x_1),
\end{align}
for some $f\in k[[x]]$.
    \end{enumerate}

\end{definition}

\begin{rmk}
The term \textit{non-commutative} here should be understood to mean \textit{non-cyclic lift}. The name \textit{potential function} is borrowed from non-commutative geometry.
\end{rmk}

In what follows, we study the relation between the (NCkrv2) and the original (krv2) equations. Let $N$ be the symmetrization map, defined to be, for a homogeneous element $s_1\ldots s_m$ of degree $m$ with $s_i\in \{x_0,x_1\}$,
\begin{align*}
N : |A| & \to A \\
|s_1\ldots s_m| & \mapsto \sum^m_{i=1}s_i\ldots s_{i-1+m},
\end{align*}
where the indices are cyclic modulo $m$. An element in the image of $N$ is said to be \emph{cyclic invariant}. Notice that $|N(|a|)|=m|a|$. 

\begin{proposition}[\cite{ANR}, Section 3]
The Fox pairing $\rho_{\rm{KKS}}$ induces a Lie bracket on the space of cyclic words $|A|$, 
\begin{equation*}
\{|a|,|b|\}_{\rm{necklace}}:=|b'S\left(\rho_{\rm{KKS}}(a'',b'')'\right)a'\rho_{\rm{KKS}}(a'',b'')''|,
\end{equation*}
and the reduced coaction induces a Lie cobracket on the space $|A|$,
\begin{align*}
\delta_{\rm{necklace}}:|A|\ & \to |A|\otimes |A| \\ 
                            |a| & \mapsto |a'S(\mu(a''))'|\otimes |\mu(a'')''|-|\mu(a'')''|\otimes |a'S(\mu(a''))'|,
\end{align*} 
where we use the Sweedler notation for the coproduct.
Moreover, they coincide with the necklace bracket and cobracket associated to the star-shaped quiver introduced by Schedler \cite{Schedler2005}.

\end{proposition}

We now relate the necklace bracket and the reduced coaction to the Lie bracket of the special derivations and the divergence map.

\begin{lemma}[\cite{GTgenus0}, Lemma 8.3; \cite{Goncharov2001}, Proposition 5.1] \label{lemma:hamiltonian}
The map 
\begin{equation} \label{eq:H-function}
    \begin{split}
        H: |A| & \to \rm{sDer}(A) \\ 
        |a| & \mapsto \left(-d^R_{0}N(|a|), -d^R_{1}N(|a|)\right)
    \end{split}
\end{equation}
is an isomorphism between the Lie algebras $(|A|/k\cdot 1$,$\{-,-\}_{\rm{necklace}})$ and $(\rm{sDer}(A),[-,-])$. The inverse map $H^{-1}$ maps a homogeneous element $u=(u_0,u_1)$ of degree $m-1$ to its Hamiltonian function $-\frac{1}{m}|x_0~u_0~+~x_1~u_1|$.
\end{lemma}

\begin{lemma}\label{lemma:div_mu}
	If $|a|\in |A|$ is homogeneous of degree $m$, then
	$
		\text{div}(H(|a|))=  -\frac{1}{m-1}| \mu(N(|a|))|,
	$ 
where $H$ is the map defined \eqref{eq:H-function}.
\end{lemma}

\begin{proof}
Let $H(|a|)=(u_0,u_1)\in \rm{sDer}(A)$.
The divergence \eqref{eq:divergence} is defined with respect to the leftmost cut of the Hamiltonian function $|x_0u_0+x_1u_1|$, and $\mu$ is defined for every adjacent letter in the Hamiltonian function. As the word $|a|$ is cyclic and homogeneous of degree $m$, the sum of all adjacent cuts is $(m-1)$ times the leftmost cut.
\end{proof}

\begin{proposition} \label{Prop:nkrv_krv}
Let $\psi \in \fr_k(x_0,x_1)$ and let $h_\psi$ be its potential function. 
If $h_\psi$ is cyclic invariant and satisfies the (NCkrv2) equation \eqref{eq:nonKV}, then $(-d^R_0(h_{\psi}),-d^R_1(h_{\psi}))\in \mathfrak{krv}_2$.
\end{proposition}

\begin{proof}
As each inhomogeneous series decomposes as the sum of its homogeneous components, it suffices to prove the statement for a homogeneous element $\psi$, in which case $h_{\psi}$ is also homogeneous. 
Since $h_{\psi}$ is cyclic invariant, $h_{\psi}=N(|a|)$ for some $a\in A$, and by the isomorphism $H$ of Lemma \ref{lemma:hamiltonian}, $(-d^R_0(h_{\psi}),-d^R_1(h_{\psi})) \in \mathfrak{sder}_2$. Since $h_\psi$ satisfies the (NCkrv2) equation, Lemma \ref{lemma:div_mu} implies that $(-d^R_0(h_{\psi}),-d^R_1(h_{\psi}))$ satisfies the (krv2) equation, after rescaling $f$ by a factor $-\frac{1}{m-1}$. Thus, $( -d^R_0(h_{\psi}),-d^R_1(h_{\psi})) \in \mathfrak{krv}_2$.
\end{proof}

\subsection{$\mathfrak{rc}_0$ and $\mathfrak{krv}_2$}
Recall that by Theorem \ref{th:poly_rc_0}, $\psi\in \mathfrak{rc}_0$ satisfies the equation 
\begin{align} \label{eq:change1}
	\mu(\psi(x_\infty,x_1))=d^R_1(\psi(x_\infty,x_1))(x_0+x_1,0)-d^R_1(\psi(x_\infty,x_1))-d^R_1(\psi(x_\infty,x_1))(x_1,0),
\end{align}
where $x_\infty:=-x_0-x_1$.
\begin{lemma}
If $\psi\in \mathfrak{rc}_0$, then
\begin{align}\label{eq:change2}
		\mu(\psi(x_\infty,x_0))=d^R_1(\psi(x_\infty,x_1))(x_0+x_1,0)-d^R_{0}(\psi(x_\infty,x_0))-d^R_1(\psi(x_\infty,x_1))(x_0,0).
	\end{align}
\end{lemma}

\begin{proof}
Let $\sigma_{0,1}$ be the algebra automorphism of $A$ that exchanges the variables $x_0$ and $x_1$, then
$\mu\circ\sigma_{0,1}=\sigma_{0,1}\circ \mu$, 
which implies \eqref{eq:change2}.
\end{proof}

\begin{proposition}\label{prop:noncommutativekrv2}
If $\psi\in \mathfrak{rc}_0$,
then $\psi$ satisfies the (NCkrv2) equation 
for $f(x):=xd^R_1(\psi(x_\infty,x_1))(x,0)$.
\end{proposition}

\begin{proof}
	By direct computation,
	\begin{align*}
		\mu(x_0\psi(x_\infty &,x_0)+x_1\psi(x_\infty,x_1))\\
		&=x_0d^R_0(\psi(x_\infty,x_0))+x_0\mu(\psi(x_\infty,x_0)) +x_1d^R_1(\psi(x_\infty,x_1))+x_1\mu(\psi(x_\infty,x_1))\\
		&=f(x_0+x_1)-f(x_0)-f(x_1),
	\end{align*}
	where the second equality follows from equations \eqref{eq:change1} and \eqref{eq:change2}.
\end{proof}

Consider the map
\begin{align*}
    L_1: \mathfrak{fr}^{>1}_k(x_0,x_1) & \to \mathfrak{tder}_2 \\
    \psi & \mapsto \bigl(
-\psi(x_\infty,x_0),
-\psi(x_\infty,x_1)
\bigr).
\end{align*}
Under the notations introduced above, \(u_\psi:=\rho(L_1(\psi))\) is the tangential derivation defined on the generators by $u_\psi(x_i)
=-[x_i,\psi(x_\infty,x_i)]$, for $i\in \{0,1\}.$

\begin{proposition}\label{prop:Hex-Lie-subalgebra}

The vector space
\[
\mathrm{Vep}
=
\left\{
\psi\in\mathfrak{fr}^{>1}_k(x_0,x_1)
\;\middle|\;
[x_1,\psi(x_\infty,x_1)]
+
[x_0,\psi(x_\infty,x_0)]
=0
\right\}
\]
is a Lie subalgebra of $\frone$ for the Ihara bracket \eqref{eq:Ihara_bracket}.
Moreover, the map
$L_1:
\mathrm{Vep}\to\mathfrak{sder}_2$ is an injective Lie algebra homomorphism, 
i.e. for all \(\psi_1,\psi_2\in\mathrm{Vep}\),
\[
L_1(\{\psi_1,\psi_2\}) = [L_1(\psi_1),L_1(\psi_2)].
\]
\end{proposition}

\begin{proof}
We start by rewriting the Vep condition. Let 
\(\psi\in\mathfrak{fr}^{>1}_k(x_0,x_1)\), then
\begin{align*}
u_\psi(x_\infty)
&=
-u_\psi(x_0)-u_\psi(x_1)=
[x_0,\psi(x_\infty,x_0)]
+[x_1,\psi(x_\infty,x_1)],
\end{align*}
hence
\[
\psi\in\mathrm{Vep}
\quad\Longleftrightarrow\quad
u_\psi(x_\infty)=0
\quad\Longleftrightarrow\quad
u_\psi\in\mathfrak{sder}_2.
\]

Let \(\psi_1,\psi_2\in \mathrm{Vep}\), then \(u_{\psi_j}(x_\infty)=0\)
for $j \in \{1,2\}$ and it follows that for 
every
\(\varphi\in\mathfrak{fr}_k(x_0,x_1)\) and
\(i\in\{0,1\}\), we have
\begin{equation}\label{eq:substitution-derivation}
u_{\psi_j}\bigl(\varphi(x_\infty,x_i)\bigr)
=
-(d_{\psi_j}\varphi)(x_\infty,x_i),
\end{equation}
where $d_{\psi_j}$ is the derivation defined in \eqref{eq:Ihara_bracket}.
Indeed, after the substitution $x_0 \mapsto x_\infty, x_1 \mapsto x_i$ (for $i\in \{0,1\}$), the two sides of
\eqref{eq:substitution-derivation} are derivations which agree on the two free
generators: the first is sent to
\(u_{\psi_j}(x_\infty)=0\), while the second is sent to
$u_{\psi_j}(x_i)
= -[x_i,\psi_j(x_\infty,x_i)].$

For \(i\in\{0,1\}\), we compute
\begin{align*}
[u_{\psi_1},u_{\psi_2}](x_i)
={}&
- u_{\psi_1}\bigl([x_i,\psi_2(x_\infty,x_i)]\bigr)
+
u_{\psi_2}\bigl([x_i,\psi_1(x_\infty,x_i)]\bigr)                 \\
={}&
\bigl[[x_i,\psi_1(x_\infty,x_i)],
      \psi_2(x_\infty,x_i)\bigr]-
\bigl[[x_i,\psi_2(x_\infty,x_i)],
      \psi_1(x_\infty,x_i)\bigr]                      +
\left[
x_i,
\bigl(
d_{\psi_1}(\psi_2)-d_{\psi_2}(\psi_1)
\bigr)(x_\infty,x_i)
\right] \\
={}& \left[
x_i,
[\psi_1,\psi_2](x_\infty,x_i)
\right] + \left[
x_i,
\bigl(
d_{\psi_1}(\psi_2)-d_{\psi_2}(\psi_1)
\bigr)(x_\infty,x_i)
\right],
\end{align*}
where the third equality follows from the Jacobi identity.
Therefore,
\begin{align*}
[u_{\psi_1},u_{\psi_2}](x_i)
&=
\left[
x_i,
\bigl(
[\psi_1,\psi_2]
+d_{\psi_1}(\psi_2)
-d_{\psi_2}(\psi_1)
\bigr)(x_\infty,x_i)
\right]                                                     \\
&=
\left[
x_i,
-\{\psi_1,\psi_2\}(x_\infty,x_i)
\right],
\end{align*}
by definition of the Ihara bracket.
Hence  $\rho([L_1(\psi_1), L_1(\psi_2)]) = \rho(L_1(\{\psi_1, \psi_2\}))$, and since $\rho$ is injective on $L_1(\fr_k^{>1})$, we conclude $[L_1(\psi_1),L_1(\psi_2)]$ 
$=L_1(\{\psi_1,\psi_2\}).$ 
Moreover, it follows from the previous computations that
\begin{align*}
    u_{\{\psi_1, \psi_2\}}(x_\infty) = [u_{\psi_1},u_{\psi_2}](x_\infty) = u_{\psi_1}(u_{\psi_2}(x_\infty)) - u_{\psi_2}(u_{\psi_1}(x_\infty)) = 0,
\end{align*}


which is equivalent to
$\{\psi_1,\psi_2\}\in\mathrm{Vep}.$ Thus, \(\mathrm{Vep}\) is a Lie subalgebra 
and $L_1$ is a Lie algebra homomorphism. Finally, it is injective because the change of variables
\[
\psi(x_0,x_1)
\mapsto
\psi(-x_0-x_1,x_1)
=
\psi(x_\infty,x_1)
\]
is an involutive automorphism of
\(\mathfrak{fr}_k(x_0,x_1)\).
\end{proof}

We conclude this section with our main result, Theorem \ref{th:Kashiwara--Vergne-E}.
\begin{theorem}[Theorem \ref{th:Kashiwara--Vergne-E}]
\label{th:Kashiwara--Vergne}
The maps $L$ and $L_1$
\begin{align*}
        L: \mathfrak{dmr}_0 \cap \mathrm{Skew} & \to \mathfrak{rc}_0\\
        \psi & \mapsto \psi,
    \end{align*}
\begin{align*}
\mathfrak{dmr}_0\cap \mathrm{Skew}\cap\mathrm{Vep}&\xrightarrow{L} \mathfrak{rc}_0\cap \mathrm{Vep}\xrightarrow{L_1} \mathfrak{krv}_2&\\
\psi(x_0,x_1)&\mapsto \psi(x_0,x_1)\mapsto (-\psi(-x_0-x_1, x_0), -\psi(-x_0-x_1, x_1)),
\end{align*}
are injective Lie algebra maps.

\end{theorem}

\begin{proof}

First, every space in the chain is a Lie algebra; the relevant results are \cite{Racinet2002} for $\dmr_0$, Lemma \ref{lem:Skew-Lie-subalgebra} for Skew, Proposition \ref{prop:Hex-Lie-subalgebra} for Vep, and Theorem \ref{theorem:property} for $\rc_0$. \\

The first map $L$ 
is an injective Lie algebra morphism by Theorem \ref{theorem:rc_dmr}, which shows that if $\psi \in \mathfrak{dmr}_0\cap$Skew then $\psi\in \mathfrak{rc}_0$. 
For the second map $L_1$, an element $\psi \in \mathfrak{rc}_0$ satisfies the (NCkrv2) equation by Proposition \ref{prop:noncommutativekrv2}. For $\psi\in \mathrm{Vep}$, this implies
\begin{equation*}
    h_{\psi}=x_1\psi(x_{\infty},x_1)+x_0\psi(x_{\infty},x_0)=\psi(x_{\infty},x_0)x_0+\psi(x_{\infty},x_1)x_1,
\end{equation*}
i.e. that $h_{\psi}$ is cyclic invariant. It then follows from Proposition \ref{Prop:nkrv_krv} that $(-d^R_0(h_{\psi}),-d^R_1(h_{\psi}))=(-\psi(x_{\infty},x_0),-\psi(x_{\infty},x_1))$ lies in $\mathfrak{krv}_2$. Finally, it is an injective Lie algebra homomorphism by Proposition \ref{prop:Hex-Lie-subalgebra}.

\end{proof}


\bibliographystyle{cas-model2-names}

\bibliography{double_shuffle}



\end{document}